%% file: prism-z.tex
\documentclass[12pt, reqno]{amsart}
\RequirePackage[all,ps,cmtip]{xy}
\usepackage{amssymb}
\usepackage[all]{xy}
\usepackage{mathrsfs}
\usepackage{enumerate}
\usepackage{amscd}
\usepackage{enumitem}
\usepackage{bbm}
\usepackage{setspace}
\usepackage{eucal}
\usepackage{xargs}
\usepackage[usenames]{color}
\usepackage{fullpage,url,mathrsfs,stmaryrd}
\usepackage{amsmath,amsthm,amssymb,mathrsfs,stmaryrd,color}
\usepackage[utf8]{inputenc}
\usepackage[T1]{fontenc}
\usepackage{relsize}
\usepackage[bbgreekl]{mathbbol}
\usepackage{amsfonts}

\newtheorem{theo}[subsubsection]{Theorem}
\newtheorem*{theo*}{Theorem}
\newtheorem*{lemm*}{Lemma}
\newtheorem{coro}[subsubsection]{Corollary}
\newtheorem{lemm}[subsubsection]{Lemma}
\newtheorem{prop}[subsubsection]{Proposition}
\theoremstyle{definition}
\newtheorem{rema}[subsubsection]{Remark}

\input{commands.tex}

\usepackage{relsize}

\title{Prismatization over $\Z$}

\begin{document}
	
\author{Lance Gurney}

\begin{abstract} The aim of this article is to given an extension of the prismatization functor for $p$-adic formal schemes (whose construction was first sketched by Drinfeld in \cite{Drinfeld} and then given by Bhatt-Lurie in \cite{Bhatt-Lurie}) to all schemes over $\Spec(\Z)$. We then prove some basic properties of this extension (algebraicity, flatness for syntomic morphisms, perfectness of cohomology) and show that for smooth schemes over $\Q$ this construction recovers (a version of) the filtered de Rham stack.
\end{abstract}

\maketitle

\tableofcontents

\section*{Introduction} Fix a prime $p$. In the remarkable articles \cite{Bhatt-Lurie} and \cite{bhatt-lurie-ii} Bhatt--Lurie introduced a certain functor called prismatization: \[\Sch_{\Z_p}\to  \St_{\Sigma_p}: X\mapsto X^{\Prism_p}\] from the category of $p$-adic formal schemes to the category of stacks over a $p$-adic formal stack $\Sigma_p:=\Spf(\Z_p)^{\Prism_p}$.\footnote{The construction of the stack $\Sigma_p$ and the functor $X\mapsto X^{\Prism_p}$ was first sketched by Drinfeld in \cite{Drinfeld} and we have opted to use his slightly more compact notation for these functors, modified with the subscript $p$ to denote its dependence on the prime $p$.}

The prismatic cohomology of a $p$-adic formal scheme $f: X\to \Spf(\Z_p)$ is given by \[H_{\Prism_p}(X):=R f^{\Prism_p}_*(\mc{O}_{X^{\Prism_p}})\in  D(\Sigma_p)\] Prismatic cohomology has many miraculous properties chief among which is that it determines essentially all other known $p$-adic cohomology theories. For example, there is a pair of natural maps $i_{\dR}, i_{HT}:\Spf(\Z_p)\to \Sigma_p$ and when $X$ is nice enough its de Rham cohomology and Hodge-Tate cohomology are recovered via pull-back along these maps: \[Li^*_{\dR}(H_{\Prism_p}(X))\isomto H_{\dR}(X)\in D(\Z_p)^{\wedge} \quad \text{ and } \quad Li^*_{HT}(H_{\Prism_p}(X))\isomto H_{HT}(X)\in D(\Z_p)^{\wedge}.\]

The aim of this article is to construct an integral version of the $p$-adic prismatization functor above. That is we construct a functor \begin{equation}\Sch_\Z\to \St_{\Sigma_\Z}: X\mto X^{\Prism_\Z}\label{eq:pris-functor}\end{equation} from the category of all schemes to the category of stacks over $\Sigma_\Z:=\Spec(\Z)^{\Prism_\Z}$ and following Bhatt-Lurie if $f: X\to \Spec(\Z)$ is a map we define the integral prismatic cohomology of $X$ to be \[H_{\Prism_\Z}(X):=Rf_*^{\Prism_\Z}(\mc{O}_{X^{\Prism_\Z}})\in D(\Sigma_\Z)\]

The main theorem of this article is that we can recover $p$-adic prismatization and $p$-adic prismatic cohomology from the integral prismatization and integral prismatic cohomology:

\begin{theo*} For $X$ quasi-compact, quasi-separated and syntomic over $\Z$ there is a natural isomorphism \[X^{\Prism_\Z}\times _{\Spec(\Z)}\Spf(\Z_p)\isomto X_{\Z_p}^{\Prism_p}\] where $X_{\Z_p}=X\times_{\Spec(\Z)}\Spf(\Z_p)$ is the $p$-adic completion of $X$. Moreover, writing $i_p: \Sigma_p=\Sigma_{\Z}\times_{\Spec(\Z)}\Spf(\Z_p)\to \Sigma_\Z$ for the induced map, there is a natural identification: \[Li^*_p(H_{\Prism_\Z}(X))\isomto H_{\Prism_p}(X_{\Z_p})\in D(\Sigma_p).\]
\end{theo*}

In words, the theorem above says that the $p$-adic completions of the integral prismatic cohomologies of $X/\Z$ are the $p$-adic prismatic cohomologies of the $p$-adic completions of $X$. A natural question is then what is the rationalisation of the integral prismatic cohomology? The answer is given by the de Rham cohomology of $X_\Q:=X\times_{\Spec(\Z)}\Spec(\Q)$ equipped with the Hodge filtration.

To understand this, let us say a little about de Rham cohomology of smooth schemes over $\Q$. Recall that if $X/\Q$ is smooth then $\Fil_H^*H_{\dR}(X)$, the Hodge filtered de Rham cohomology of $X$, is naturally an object of $DF(\Q)^{\wedge}$, the filtered complete category of complexes of $\Q$-vector spaces. The Rees equivalence \[\mathrm{Rees}: D((\widehat{\A}^1/\Gm)_{\Q})\isomto DF(\Q)^{\wedge}\] allows us to identify $DF(\Q)^{\wedge}$ with $D((\widehat{\A}^1/\Gm)_{\Q})$, the derived category of quasi-coherent sheaves on $(\widehat{\A}^1/\Gm)_\Q=\widehat{\A}^1/\Gm\times_{\Spec(\Z)}\Spec(\Q)$. Therefore, we may view $\Fil_H^*H_\dR(X)$ as an object of $D((\widehat{\A}^1/\Gm)_\Q)$.

Along with the prismatization functor, we introduce a `de Rhamification' functor\footnote{Other `de Rhamification' functors have been considered previously, notably by Simpson \cite{simpson}. Ours is similar to his, although it differs in some important regards; see (\ref{rema:other-approaches}) for details on this.} \begin{equation}\Sch_\Q\to \St_{(\widehat{\A}^1/\Gm)_\Q}: X\mapsto X^{\dR}\label{eq:de-rhamification}\end{equation} satisfying $\Spec(\Q)^{\dR}=(\widehat{\A}^1/\Gm)_\Q$, and show that if $f: X\to \Spec(\Q)$ is smooth, quasi-compact and quasi-separated then \[Rf_*^{\dR}(\mc{O}_{X^\dR})\isomto \Fil_H^*H_\dR(X)\in D((\widehat{\A}^1/\Gm)_\Q).\]

The relation between rational prismatic cohomology and de Rham cohomology is given by the following:

\begin{theo*} For $X$ quasi-compact, quasi-separated and smooth over $\Z$, we have \[X^{\Prism_\Z}\times_{\Spec(\Z)}\Spec(\Q)\isomto X^{\dR}_\Q\] where $X_\Q=X\times_{\Spec(\Z)}\Spec(\Q)$. Moreover, writing $i_\Q: (\widehat{\A}^1/\Gm)_\Q \to \Sigma_\Z$ for the induced map, there is a natural identification: \[Li_\Q^*(H_{\Prism_\Z})\isomto \Fil_H^*H_{\dR}(X_\Q)\in D((\widehat{\A}^1/\Gm)_\Q).\]
\end{theo*}

Let us explain the outline of the article.

In Section \ref{sec:de-rham} we give the construction of the de Rhamification functor (\ref{eq:de-rhamification}) and show that for smooth schemes over $\Q$ the quasi-coherent cohomology of $X^{\dR}$ coincides with the Hodge filtered de Rham cohomology.

In Section \ref{sec:witt-vec} we quickly review some basic properties of the big Witt vectors $W$ and construct a certain invertible $W$-module $\mc{V}$.

In Section \ref{sec:sigma} we introduce the stack $\Sigma_\Z$ mentioned earlier in the introduction. Our construction of $\Sigma_\Z$ directly parallels the $p$-adic theory where $\Sigma_p$ is defined as the quotient of a $p$-adic formal subscheme $W_{\dist}^{(p)}\subset W^{(p)}$ of the $p$-typical Witt vectors $W^{(p)}$ by the action of the group of units $W^{(p), \times}$. Thus, the key definition in this section is of a certain formal subscheme $W_\dist\subset W$. Once this is done we define $\Sigma_\Z:=W_\dist/W^\times$, prove some algebraicity properties and show that that its $p$-adic completions coincide with the stacks $\Sigma_p$.

Finally, in Section \ref{sec:prismatization} we introduce the prismatization functor $X\mapsto X^{\Prism_\Z}$. As with $\Sigma_\Z$, our construction of this functor directly parallels that of the $p$-adic case. One the prismatization functor is defined, we prove various algebraicity and flatness properties. Finally, we prove the two theorems above.

\section{Preliminaries}

\subsection{Conventions}

\subsubsection{Stacks} We denote by $\aS$ the $\infty$-category of animated sets (equivalently the $\infty$-category of spaces). We denote by $\Aff$ the opposite of the 1-category of rings $\Ring$ and by $\dAff$ the opposite of the $\infty$-category of animated rings $\aRing$. For a short overview of animation and related concepts see \S 5.1 of \cite{scholze-ces-flat}.

We denote by $\dPSt$ the category of derived prestacks, i.e. the category of functors $X: \dAff^\circ\to \aS$. Equipping $\dAff$ with the fpqc topology, we have the full subcategory $\dSt\subset \dPSt$ of derived stacks.

Inside $\dSt$ we have the full subcategory of classical stacks $\St$ (which we will just call stacks) which is equivalent to the category of functors $\Aff^\circ\to \aS$ satisfying the sheaf condition for fpqc covers; equivalently they are those derived stacks which can be written as a colimit of non-derived affine schemes. Note that the inclusion $\St\to \dSt$ commutes with colimits but not with limits (unless they are `Tor-independent').

We will not really have any use for the rather wild generality above, but it is convenient to consider all of the objects which will appear: schemes, derived affine schemes, 1-stacks and so on, as living (fully faithfully) within a single larger category.

All limits will mean limits taken in $\St$ and not their derived counterparts, unless explicitly stated.

\subsubsection{Algebraic stacks}\label{subsubsec:algebraic-stacks} A morphism of $f: \mc{X}\to \mc{Y}$ of stacks is called schematic if for all schemes $Y$ and all morphisms $Y\to \mc{Y}$ the fibre product $\mc{X}\times_{\mc{Y}}$ is representable by a scheme. If $f$ is schematic and $P$ is any property of morphisms of schemes stable under base change then we say that $f$ has property $P$ if the morphisms $f\times_{\mc{Y}}Y$ have property $P$ for all maps $Y\to \mc{Y}$ with $Y$ a scheme.\footnote{Foundationally speaking, schematic morphisms are not so useful; for instance the property of being schematic does not descend. However, this notion will suffice for our purposes.}

Schematic morphisms which tend to arise in nature are affine morphisms and open immersions (and compositions thereof).

A stack $\mc{X}$ is called algebraic if there exists a schematic morphism $X\to \mc{X}$ which is representable by faithfully flat morphisms of schemes which are covering maps for the fpqc topology; we call $X\to \mc{X}$ a presentation of $\mc{X}$.

A morphism of stacks $\mc{X}\to \mc{Y}$ is called algebraic if for all schemes $Y\to \mc{Y}$, the stack $\mc{X}\times_{\mc{Y}} Y$ is algebraic. Any morphism of algebraic $\mc{Y}$-stacks is algebraic and algebraic morphisms are stable under composition.

An algebraic morphism of stacks $f: \mc{X}\to \mc{Y}$ is called (faithfully) flat if for all maps $Y\to \mc{Y}$ with $Y$ a scheme, and all presentations $X\to \mc{X}\times_{\mc{Y}}Y$ the composition $X\to \mc{X}\times_{\mc{Y}}Y\to Y$ is a flat morphism of schemes (and a covering map for the fpqc topology). Note that (faithfully) flat algebraic morphisms are stable under composition.

Finally, we note that if $\mc{X}\to \mc{Y}$ is a flat algebraic morphism and $\mc{Z}\to \mc{Y}$ is any morphism the the limit $\mc{X}\times_{\mc{Y}}\mc{Z}$ is Tor-independent i.e. it agrees with the limit computed in $\dSt$.

\subsubsection{Ring schemes and invertible modules} If $S$ is a stack and $\mc{R}$ is a sheaf of rings over $S$ then we write $\mc{R}^\times$ for the sheaf of units of $\mc{R}$. If $S'$ is an $S$-scheme then an $\mc{R}_{S'}$-module $\mc{P}$ is called invertible if it is fpqc locally isomorphic to $\mc{R}_{S'}$. We identify the $S'$-points of the classifying $S$-stack $B \mc{R}^\times$ with the groupoid of invertible $\mc{R}_{S'}$-modules and we identify the $S'$-points of the $S$-stack $\mc{R}/\mc{R}^\times$ with the groupoid whose objects are invertible $\mc{R}_{S'}$-modules $\mc{P}$ equipped with an $\mc{R}_{S'}$-linear map $\xi: \mc{P}\to \mc{R}_{S'}$.

\subsubsection{Quasi-ideals, cones and ring stacks} Here we give a brief over view of the notion of the cone construction and the notion of ring stacks. For a thorough discussion of this concept we refer the reader to \cite{Drinfeld-rings}.

If $R$ is an $A$-algebra, a quasi-ideal in $R$ is an $R$-module $I$ equipped with an $R$-linear map $d: I\to R$ such that $x d(y)=y d(x)$ for all $x, y\in I.$ We denote by $\Cone(d:I\to R)$ or just $\Cone(d)$ the groupoid whose set of object is $R$ and whose morphisms are given by $\Hom(r_1, r_2):=d^{-1}(r_2-r_1)$ for $r_1, r_2\in R$. Equivalently, we may view $I$ as acting on the set $R$ via translation through $d$, in which case $\Cone(d)$ is equivalent to the quotient groupoid of $R$ by this action of $I$.

The groupoid $\Cone(d)$ inherits the structure of an animated $A$-algebra in $\aS$ as follows. Denote by $R_n$ the $n$-fold fibre product of the map $R\to \Cone(i)$ with itself. Each $R_n$ is a discrete groupoid which can be explicitly described as the set $R\times I^{n-1}$. Equipping this set with the obvious addition law and with the multiplication given by \[(r, x_1, \ldots, x_{n-1})\cdot (s, y_1, \ldots, y_{n-1})=(rs, ry_1+sx_1+d(x_1)y_1, \ldots, ry_{n-1}+sx_{n-1}+d(x_{n-1})y_{n-1})\] makes $R_n$ into a commutative ring, which is also an $A$-algebra for the obvious action of $A$.

The natural maps $R_n\to R_m$ for $n, m\in \N$ are $A$-linear homomorphisms and we obtain a simplicial $A$-algebra $R_\bullet$. Taking the colimit of this diagram in $\aRing_A$ yields an animated $A$-algebra whose underlying groupoid is equivalent to the groupoid $\Cone(d: I\to R)$. Note that if $d$ is injective then $\Cone(d)\isomto R/I$.

More generally, of $\mc{R}$ is a sheaf of $A$-algebras over some stack $S$ and $i: \mc{I}\to \mc{R}$ is a sheaf if quasi-ideals, then we can form the quotient stack $\Cone(d):=\Cone(d:\mc{I}\to \mc{R})$. Sheafifying the construction above we may view the stack $\Cone(d)$ as representing a functor: \[\Cone(d): \Aff^\circ_S \to \aRing_A.\]

We call such an object an $A$-algebra stack over $S$.

An example which will occur frequently is the following: if $\mc{R}$ is a sheaf of rings over a stack $S$ and $\mc{I}$ is an invertible $\mc{R}$-module then any $\mathcal{R}$-linear map $d: \mc{I}\to \mc{R}$ defines a quasi-ideal in and hence we have the associated ring stack $\Cone(d)$ over $S$.

\subsubsection{Derived categories of quasi-coherent sheaves}\label{subsec:quasi-coherence} For a (possibly) derived stack we denote by $D(X)\subset D(X, \mc{O}_X)$ the derived $\infty$-categories of quasi-coherent sheaves and sheaves of $\mc{O}_X$-modules on the big derived affine site $\dAff_X$ of $X$.

For any morphism of stacks $f: X\to Y$ there is an induced pullback functor $Lf^*: D(Y, \mc{O}_Y)\to D(X, \mc{O}_X)$ whose right adjoint we denote by $Rf_+$. The functor $Lf^*$ restricts to a functor $Lf^*:D(Y)\to D(X)$ and we denote its right adjoint by $Rf_*$. While the formation of $Rf_+$ commutes with arbitrary base change of derived stacks this is not true for $Rf_*$. However, for all $\mc{M}\in D(X)$ there is a natural map $Rf_*(\mc{M})\to Rf_+(\mc{M})$ which is an isomorphism whenever $Rf_+(\mc{M})$ is quasi-coherent. There are two important cases when $Rf_+(\mc{O}_X)$ is quasi-coherent:
\begin{enumerate}[label=\textup{(\roman*)}]
\item If $X$ can be written as a finite colimit of stacks $X_i\to X$ such that each $Rf_{i+}(\mc{M}_i)$ is quasi-coherent, where $f_i=f|_{X_i}$ and $\mc{M}$ is the pullback of $\mc{M}_i$ to $X_i$, then so is $Rf_+(\mc{M})$.\footnote{If $X$ is a stack and $X_i\subset X$ for $i\in I$ is a finite collection of substacks which cover $X$ then $X$ is isomorphic to the colimit over the finite category of finite fibre products of $X_i$ for $i\in I$.}

\item If $X\to Y$ can be written as a simplicial colimit of relatively flat and affine $Y$-stacks then $Rf_+(\mc{O}_X)$ is quasi-coherent.
 \end{enumerate}
\section{de Rham cohomology over $\Q$}\label{sec:de-rham} The aim of this section is to give a stacky approach to the de Rham cohomology of smooth schemes over $\Q$. We say `a' stacky approach as there is another slighter different approach due to Simpson \cite{simpson} (see (\ref{rema:other-approaches}) for details on the differences).

\subsection{Reminders on filtrations and the Cartier stack}
\label{subsubsection:cartier-divisors} We recall that the Cartier stack $\A^1/\Gm$ is the stack over $\Spec(\Z)$ whose value on a ring $R$ is given by the groupoid whose objects are $(\eta: L\to R)$ where $L$ is a projective rank one $R$-module and $\eta$ is an $R$-linear map. We denote by $\mc{O}(1)$ the tautological line bundle on $\A^1/\Gm$ whose value on $(\eta: L\to R)$ is $L^\vee$ and by $\mc{O}(i):=\mc{O}(1)^{\otimes n}$ for $n\in \Z$. Thus, there is a canonical map \[t: \mc{O}(-1)\to \mc{O}\] of line bundles over $\A^1/\Gm$. For $\mc{M}\in D(\A^1/\Gm)$ we write $\mc{M}(i)=\mc{M}\otimes \mc{O}(i)$.

It is an observation due to Simpson \cite{simpson} that the Cartier stack is intimately related to filtrations. In current language, this is expressed via the Rees equivalence which, for each ring $R$ is an equivalence of $\infty$-categories: \[\mathrm{Rees}: D((\A^1/\Gm)\times \Spec(R))\isomto DF(R).\] The right hand side denotes the filtered derived $\infty$-category of the ring $R$, $DF(R):=\Fun(\Z^{\circ}, D(R))$ where $\Z^\circ$ denotes the ordered set $\Z$. This equivalence is symmetric monoidal if the target is given the symmetric monoidal structure given by Day convolution.\footnote{The Rees functor here is normalised so that $\Rees(\mc{O}(-1))$ corresponds to the functor $F: \Z^\circ\to D(R)$ whose value is $F(i)=R$ and $i\leq 1$ and $F(i)=0$ for $i>1$.} (For proofs of this and the following statements statements we refer the reader to \cite{HKR} and \cite{Mou}).

If $M: i\mapsto M(i)$ is an object of $DF(R)$ then we denote by $M\{n\}$ the object obtained by shifting the filtration to the left by $n$: $M\{n\}: i\mapsto M(i+n).$ Under the Rees equivalence this corresponds to tensoring with $\mc{O}(-n):=\mc{O}(1)^{\otimes -n}$.

An object $M: i\mapsto M(i)$ in $DF(R)$ is said to be complete if $R\lim_i M(i)=0$ and we write $DF(R)^{\wedge}\subset DF(R)$ for the full subcategory of complete objects. The Rees equivalence restricts to an equivalence \[\Rees: D((\widehat{\A}^1/\Gm)\times \Spec(R))\isomto DF(R)^{\wedge}.\] We call $\widehat{\A}^1/\Gm$ the formal Cartier stack.

\subsection{The de Rham stack} In this and the following section we will work solely over $\Q$ and will omit it from our notation: e.g. $\A^1$ will denote $\A^1_\Q$ etc.

Over the formal Cartier stack $\widehat{\A}^1/\Gm$ we have the smooth, affine algebraic group $\Ga(-1)$ corresponding to the line bundle $\mc{O}(-1)$ and the map \[t: \Ga(-1)\to \Ga\] corresponding to $t: \mc{O}(-1)\to \mc{O}$.

Viewing $\Ga$ as a ring scheme and $\Ga(-1)$ as a $\Ga$-module scheme, the map $t$ above defines a quasi-ideal and we denote by $\Ga^{\dR}$ the $\Q$-algebra stack over $\widehat{\A}^1/\Gm$ given by \[\Ga^{\dR}=\Cone(\Ga(-1)\to \Ga).\]

\begin{lemm}\label{lemm:props-of-gadr} Let $R$ be a ring and let $(\eta: L\to R) \in \A^1/\Gm(R)$.
\begin{enumerate}[label=\textup{(\roman*)}]
\item The natural map $\Cone(\eta: L\to R) \isomto \Ga^{\dR}(\eta: L\to R)$ is an isomorphism.
\item If $R\to R'$ is a homomorphism then \[\Ga^{\dR}(R)\otimes_R^{\mathbf{L}} R'\isomto \Ga^{\dR}(R').\]
\item The functor on categories of animated \'etale algebras: \[\mathrm{Et}_{R}\isomto \mathrm{Er}_{\Ga^{\dR}(R)}: R'\mapsto \Ga^{\dR}(R)\otimes_{R}^{\mathbf{L}}R'\] is an equivalence.
\end{enumerate}
\end{lemm}
\begin{proof}
(i) This follows from the fact that $H^1(\Spec(R), \Ga(-1))=0$.

(ii) This follows from (i).

(iii) By (i) the map of animated rings $R\to \Ga^{\dR}(R)$ is surjective on $\pi_0$ with nilpotent kernel, which gives the claim by the invariance of categories of animated \'etale algebras under such maps.
\end{proof}

\subsubsection{}\label{subsec:de-rham-stacks} If $X$ is a scheme over $\Q$ we define the de Rham stack $X^{\dR}$ of $X$ to be the prestack over $\widehat{\A}^1/\Gm$ whose value on $\Spec(R)\to \widehat{\A}^1/\Gm$ is \[X^{\dR}(\Spec(R))=X(\Ga^{\dR}(R)).\] It is clear from this definition that the functor $X\mto X^{\dR}$ commutes with Tor-independent limits.

\subsubsection{}\label{subsec:dr-fibres} Let us briefly explain what the fibres of a morphism $\Spec(B)^{\dR}\to \Spec(A)^{\dR}$ look like. If $\Spec(R)\to \Spec(A)^{\dR}$ is a map corresponding to $h: A\to \Ga^{\dR}(R)$ then the fibre product $\Spec(B)^{\dR}\times_{\Spec(A)^{\dR}}\Spec(R)$ is the $\Spec(R)$-stack whose value on an $R$-algebra $R'$ is given by \[(\Spec(B)^{\dR}\times_{\Spec(A)^{\dR}}\Spec(R))(R')\isomto \Hom_A(B, \Ga^{\dR}(R'))\] where $\Ga^{\dR}(R')$ has the $A$-algebra structure coming from $h$.

Before we consider the general properties of de Rham stacks we leave the following as an exercise for the reader:

\begin{lemm}\label{lemm:schematic} Let $S$ be a scheme, $T$ be a stack and let $X\mto \wt{X}$ be a functor from $\Sch_S\to \St_T$ which sends affine open immersions (resp. covers) to affine open immersion (resp. covers) and commutes with base change along such maps. If $P$ is a property of morphisms of schemes such that
\begin{enumerate}[label=\textup{(\alph*)}]
\item $P$ is Zariski local on the source and target, and
\item $X\mto \wt{X}$ sends morphisms of affine schemes with property $P$ to affine morphisms with property $P$.
\end{enumerate}
Then $X\mto \wt{X}$ sends morphisms with property $P$ to schematic morphisms with property $P.$
\end{lemm}

\begin{lemm} The functor $X\mapsto X^{\dR}$ send \'etale morphisms and covers to schematic \'etale morphisms and covers (and likewise for open immersions and covers). Moreover, $X^{\dR}$ is an \textup{fpqc} algebraic stack.
\end{lemm}
\begin{proof} The case of open immersions (which are \'etale monomorphisms) and open covers follows from the case of \'etale morphisms and covers by the compatibility with flat base change.

Thus by (\ref{lemm:schematic}) we are reduced to the case where $f: X=\Spec(B) \to Y=\Spec(A)$ is an \'etale map (resp. cover) of affine schemes.

The claim is then equivalent to the following (cf. (\ref{subsec:dr-fibres})): for all morphisms $\Spec(R)\to \Spec(A)^{\dR}$ corresponding to $h: A\to \Ga^{\dR}(R)$, the functor on $R$-algebras given by \[R'\mapsto \Hom_A(B, \Ga^{\dR}(R'))\] is representable by an \'etale $R$-algebra. However, it is clear that this functor is representable by the \'etale $R$-algebra $R_{B/A}$ corresponding to the \'etale $\Ga^{\dR}(R)$-algebra $\Ga^{\dR}(R)\otimes_A^{\mathbf{L}} B$ under the equivalence in (ii) of (\ref{lemm:props-of-gadr}) and this proves the claim for \'etale maps and covers.

It follows from (ii) of (\ref{lemm:props-of-gadr}) that $\Ga^{\dR}$ sends the Cech nerve of a faithfully flat homomorphism $R\to R'$ to the Cech nerve of the faithfully flat homomorphism $\Ga^{\dR}(R)\to \Ga^{\dR}(R')$, which implies that $X^{\dR}$ is a stack for the fpqc topology.
\end{proof}

\subsubsection{} For $\Spec(R)\to \wh{\A}^1/\Gm$ the natural homomorphism $R\to \Ga^{\dR}(R)$ induces for each scheme $X$, a map \[\pi_{X^\dR}: X\times \wh{\A}^1/\Gm\to X^{\dR}.\]

\begin{prop}\label{prop:de-rham-etale} Let $X$ be a smooth scheme over $\Q$.
\begin{enumerate}[label=\textup{(\roman*)}] 
\item If $U\to X$ is \'etale the diagram \[\xymatrix{U\times \wh{\A}^1/\Gm\ar[r]\ar[d] &  U^{\dR}\ar[d]\\
X\times \wh{\A}^1/\Gm\ar[r] & X^{\dR}}\] is cartesian.
\item $\pi_{X^\dR}$ is affine and faithfully flat
\item $X^{\dR}\to \wh{\A}^1/\Gm$ is algebraic and faithfully flat.\footnote{We will see later that $X^{\dR}\to \wh{\A}^1/\Gm$ is always algebraic.}
\end{enumerate}
\end{prop}
\begin{proof} (i) Unwinding the definitions, the claim is equivalent to the statement that for all maps $\Spec(R)\to \wh{\A}^1/\Gm$ the map \[\Hom_{X}(\Spec(R), U)\to \Hom_{X}(\Spec(\Ga^{\dR}(R)), U)\] is an equivalence. However, as $U\to X$ is \'etale this follows immediately upon noting that the map of animated rings $R\to \Ga^{\dR}(R)$ is surjective on $\pi_0$ with nilpotent kernel.

(ii) Repeated applications of (i) reduce us to the case where $X=\Ga^n$ in which case it is clear as $\pi_{(\Ga^n)^{\dR}}$ is given by the natural map  \[\Ga^n\to \Cone(\Ga(-1)\to \Ga)^n\isomto (\Ga^n)^{\dR}.\]

(iii) This follows from (ii) and the fact that $X\times \wh{\A}^1/\Gm\to \wh{\A}^1/\Gm$ is schematic.
\end{proof}

\subsubsection{A presentation of $X^{\dR}$}\label{subsec:presentation-of-xdr} If $X/\Q$ is smooth then the morphism $\pi_{X^{\dR}}$ can be used to give a simplicial presentation of the stack $X^{\dR}$ via: \[X_n^{\dR}:=(X\times \wh{\A}^1/\Gm)\times_{X^\dR} \cdots \times_{X^\dR}(X\times \wh{\A}^1/\Gm).\]

Let us describe $X_n^{\dR}$ in the case when $X=\Spec(A)$ is affine. Unwinding the definitions we see that if $\Spec(R)\to \wh{\A}^1/\Gm(R)$ corresponds to $(\eta: L\to R)$ then $X_n^{\dR}(R)$ is equivalent to the groupoid of commutative diagrams \[\xymatrix{A^{\otimes_\Q n}\ar[r]\ar[d] & \ar[d]R \\
A\ar[r] & \Ga^{\dR}(R)}\] where the vertical maps are the natural ones. Taking the fibres of the vertical maps shows that this groupoid is equivalent to the discrete groupoid of commutative squares \[\xymatrix{I_n \ar[r]\ar[d] & L\ar[d] \\
A^{\otimes_\Q n}\ar[r] & R}\] where $I_n=\ker(A^{\otimes_\Q n}\to A)$.

If there exists trivialisation $L\isomto R$ giving $\eta=t_R\in R$, then using the fact that $I_n$ is regular we find that the set of diagrams above is equivalent to the set of $t$-linear homomorphisms \[A^{\otimes_\Q n}[t, t^{-1} I_n]\to R\] where $t$ acts on $R$ via $\eta=t_R$.

Letting $\Gm$ act on $\Spec(A^{\otimes_\Q n}[t, t^{-1} I_n])\to \Spec(\Q[t])=\A^1$ by giving $t$ degree $1$ defines an affine morphism: \[\Spec(A^{\otimes_\Q n}[t, t^{-1} I_n])/\Gm \to \A^1/\Gm\] and it follows by construction that \[\Spec(A)^{\dR}_n\isomto \Spec(A^{\otimes_\Q n}[t, t^{-1} I_n])/\Gm\times_{\A^1/\Gm}{\widehat{\A}^1/\Gm}.\]

Now, the morphism \[\Spec(A^{\otimes_\Q n}[t, t^{-1} I_n])/\Gm \to \A^1/\Gm\] is flat and affine and corresponds, under the Rees equivalence, to the filtered $\Q$-algebra $\Fil_{I_n} A^{\otimes_n \Q}$ where $\Fil_{I_n}$ denotes the $I_n$-adic filtration. Therefore, $\Spec(A)^{\dR}_n$ corresponds to the filtered completion of $\Fil_{I_n} A^{\otimes_n \Q}$ which, writing $(-)^\wedge$ for $I_n$-adic completion, is given by $\Fil_{I_n^{\wedge}} (A^{\otimes_\Q n})^{\wedge}$.

\subsection{Cohomology of $X^{\dR}$}

\subsubsection{de Rham cohomology} Let $f: X\to \Spec(\Q)$ be a smooth $\Q$-scheme. We view the de Rham complex $\Omega_{X/\Q}^\bullet$ of $X$ as a complex of sheaves of $\Q$-vector spaces on the small \'etale site of $X$. Equipping $\Omega_{X/\Q}^\bullet$ with the stupid filtration defines a filtered complex of sheaves of $\Q$-vector spaces $\Fil_H^\ast \Omega^\bullet_{X/\Q}$ on the small \'etale site of $X$ and the Hodge filtered de Rham cohomology $\Fil_H^\ast H_{\dR}(X/\Q)\in DF(\Q)$ of $X$ is defined the be the hypercohomology of the filtered complex $\Fil_H^\ast \Omega^\bullet_{X/\Q}$.

\begin{prop} If $f: X\to \Spec(\Q)$ is smooth, quasi-compact and quasi-separated then $Rf_+^{\dR}(\mc{O}_{X_\dR})$ is quasi-coherent and corresponds under the Rees equivalence $D(\wh{\A}^1/\Gm)\isomto DF(\Q)^{\wedge}$ to $\Fil_{H}^*H_{\dR}(X/\Q)$.
\end{prop}
\begin{proof} The quasi-coherence of $Rf_+(\mc{O}_{X^\dR})$ follows from (\ref{subsec:quasi-coherence}) combined with the compatibility of $X\mapsto X^{\dR}$ with open covers.

For the second claim it will be enough to construct isomorphisms \[\mathrm{Rees}(Rf_+^{\dR}(\mc{O}_{X^{\dR}}))\isomto \Fil_{H}^*H_{\dR}(X/\Q)\] for $X=\Spec(A)$ affine which are functorial in \'etale morphisms. In this case, \[\Fil_{H}^*H_{\dR}(X/\Q)\isomto \mathrm{Rees}(\Omega^{\bullet}_A(\bullet))\] where $\Omega^{\bullet}_A(\bullet)\in D(\wh{\A}^1/\Gm)$ denotes the complex of flat quasi-coherent modules on $\widehat{\A}^1/\Gm$ whose terms are the flat quasi-coherent modules $\Omega^{i}_A(i)$ and whose differential is $d\otimes t^{\vee}(i):\Omega^{i}_A(i)\to \Omega^{i}_A(i+1)$. Thus it is enough to construct an isomorphism \[Rf_+^{\dR}(\mc{O}_{X^{\dR}})\isomto \Omega^{\bullet}_A(\bullet).\]

For $i\geq 0$, $\Omega_{X^\dR}^i:=R\pi_{X^{\dR}+}(\Omega_{X/\Q}^i\boxtimes \mc{O}(i))$ is a flat quasi-coherent sheaf on $X_{\dR}$ as $\pi_{X^{\dR}}$ is flat and affine. Under the Rees equivalence the flat quasi-coherent module $L\pi_{X^\dR}^*\Omega_{X^\dR}^i$ corresponds (cf. (\ref{subsec:presentation-of-xdr})) to the flat filtered complete $A$-module $\Fil_{I^\wedge}\Omega^i_A\otimes_{A} (A\otimes_\Q A)^{\wedge}\{-i\}$ where the hat denotes $I=\ker(A\otimes_\Q A\to A)$-adic completion.

The linearisation of the de Rham differential $d: \Omega_A^i\to \Omega_A^{i+1}$: \[\wt{d}:\Omega_{A}^i\otimes_\Q A=\Omega_{A}^i\otimes_A (A\otimes_\Q A) \to \Omega_A^{i+1}\] satisfies  $\wt{d}(I^j\cdot \Omega_A^i\otimes_\Q A)=0$ for $j>1$ an so it extends to a filtered map \[\Fil_{I^\wedge}\Omega_A^i\otimes_A (A\otimes_\Q A)^{\wedge}\{-i\}\to \Fil_{\mathrm{triv}}\Omega_A^{i+1}\{-i-1\}.\] The right hand side corresponds under the Rees equivalence to $\Omega_{X/\Q}^{i+1}\boxtimes \mc{O}(i+1)$ and so we have defined a map: \[L\pi_{X^\dR}^*(\Omega_{X^\dR}^i)\to \Omega_{X/\Q}^{i+1}\boxtimes \mc{O}(i+1)\] and hence a map \[d_{X^{\dR}}: \Omega_{X^\dR}^i\to \Omega_{X^\dR}^{i+1}.\]

These maps define a complex of flat quasi-coherent sheaves $(\Omega_{X^\dR}^\bullet, d_{X^\dR})$ on $X^{\dR}$ and the canonical map $\mc{O}_{X_\dR}\to \Omega_{X^\dR}^0$ induces a quasi-isomorphism \[\mc{O}_{X_\dR}\isomto \Omega_{X^\dR}^\bullet.\] This can be checked after pull-back along $\pi_{X^\dR}$ in which case it is the filtered Poincar\'e lemma (see \cite{berthelot-ogus}).

Writing $p: X\times \wh{\A}^1/\Gm\to \widehat{\A}^1/\Gm$ for the projection (which is flat and affine) we see that \[Rf_+^\dR(\Omega_{X_\dR}^i)\isomto Rp_+(\Omega_{X/Q}^i\boxtimes \mc{O}(i)) \isomto \Omega_{A}^i(i)\] while the direct imeage of the differential $d_{X^{\dR}}$ is \[d\otimes t^\vee(i): \Omega_{A}^i(i)\to \Omega_{A}^{i+1}(i+1).\] Therefore, \[Rf_{+}^{\dR}(\mc{O}_{X^\dR})\isomto \Omega_{A}^{\bullet}(\bullet).\]
\end{proof}

\begin{rema}\label{rema:other-approaches} There is another `stacky' approach to the de Rham cohomology for schemes over $\Q$ due to Simpson \cite{simpson}. In Simpon's set up one works over the (uncompleted) Cartier stack $\A^1/\Gm$ and instead of working with the ring stack $\Cone(\Ga(-1)\to \Ga)$, one works with the ring stack $\Ga^{\dR, \mathrm{Sim}}:=\Cone(\widehat{\mathbf{G}}_{\mathrm{a}}(-1)\to \Ga)$ and, using the same formulas as in (\ref{subsec:de-rham-stacks}), defines stacks $X^{\dR, \mathrm{Sim}}$. The major upshot of Simpson's approach is that one can recover unfiltered de Rham cohomology via base change along the open point $\Spec(\Q)\to \A^1/\Gm$ which is missing in the story above. However, it comes with some drawbacks:

First, the stacks $X^{\dR, \mathrm{Sim}}$ are not algebraic over $\A^1/\Gm$ due to the appearance of the formal group scheme $\wh{\mathbf{G}}_{\mathrm{a}}(-1)$ above and so the quasi-coherence of the direct images, while still true, must be proven directly.

Secondly, there is a general method (due to Mondal \cite{mondal}) by which one can try to recover stacky approaches to cohomology theories from purely cohomological information. This method can be applied in the setting outlined in this article to show that Hodge filtered de Rham cohomology of smooth affine schemes can be used to recover the stacks $X^{\dR}$, while this method fails in Simpson's approach (and will fail for any approach which recovers unfiltered de Rham cohomology).

Lastly, the theory given in this article can be extended to work for all smooth schemes over $\Z$, while it is not clear that there is a sensible way to do this in Simpson's approach (although we note that something probably can be done in this direction, see Appendice \cite{Berthelot}.)
\end{rema}
\section{Witt vectors}\label{sec:witt-vec}

\subsection{Reminders} The aim of this section is prove a handful of results regarding certain ring and ideal schemes associated to the Witt vectors. We begin with some reminders and notation (which is largely unavoidable when working with the Witt vectors).

We then construct an invertible $W$-module $\mc{V}$ where $W$ denotes the ring scheme of big Witt vectors. This invertible module comes with a canonical map $V(1): \mc{V}\to W$ which is the big Witt vector analogue of the maps $V_p(1): W^{(p)}\to W^{(p)}$ where $W^{(p)}$ denotes the $p$-typical Witt vectors and $V_p$ is the Verschiebung.

We then study the sheaf of rings $\wt{W}=W/W[F]$ obtained from $W$ by quotienting by the ideal sheaf $W[F]$ given by the intersection over all primes $p$ of the kernels of the Frobenius endomorphisms. We show that if $VW\to W$ denotes the Verschiebung ideal, then there is a natural isomorphism $VW\isomto \mc{V}\otimes_W \wt{W}$ from which we deduce various algebraicity properties of the sheaf of rings $\wt{W}$ (notably the fact that it is flat, affine and formally smooth).

For a detailed discussions of various incarnations of Witt vectors and their properties, we refer the reader to \cite{borger-witt-i}.

\subsubsection{Big Witt vectors}\label{subsec:e-typical-witt-vectors} The various versions of the Witt vectors are naturally indexed by subsets $E\subset \N$ which are closed under divisors and products of relatively prime integers; let us call these index sets. If $E$ is a finite index set then it is of the form $\{d\in \N: d|n\}$ for some $n\in \N$ and if $E$ is infinite then it is an increasing union of finite index sets.

Let us introduce some operations on index sets. If $p$ is a prime then we write $E|p=\{n/\gcd(n, p):n\in E\}\subset E$, $E_{(p)}=\{n\in E: \gcd(n, p)=1\}\subset E$, and extend this notation from the primes $p$ to all integers by iterating over prime power decompositions (noting that the order does not matter). We also write $E^{(p)}=\{1, p, p^2, \ldots\}\cap E$.

Given such a finite index set $E$, we denote by $W_E$ the flat affine ring scheme of $E$-typical Witt vectors: if $E=\{d: d|n\}$ for some $n\in \N$ then $W_E:=W_{p_1^{r_1}}\circ \cdots \circ W_{p_m^{r_m}}$ where $n=p_1^{r_1}\cdots p_m^{r_m}$ is a decomposition of $n$ into prime powers and $W_{p^r}$ denotes the $p$-typical Witt vectors of length $r+1$. If $E\subset E'$ is an inclusion of finite index sets then there are natural restriction maps $W_{E'}\to W_E$ and for $E$ infinite we set $W_E=\lim_{E'\subset E \text{ finite}} W_{E'}.$

If $E$ and $E'$ are two index sets whose elements are relatively prime then $E\cdot E'$ is again an index set and there is a natural isomorphism \[W_{E\cdot E'}\isomto W_{E}\circ W_{E'}\isomto W_{E'}\circ W_E.\] In particular, if $p\in E$ is a prime then $E^{(p)}\cdot E_{(p)}=E$ and so \[W_E\isomto W_{E^{(p)}}\circ W_{E_{(p)}}\isomto W_{E_{(p)}}\circ W_{E^{(p)}}.\]

We will mainly be concerned with the following index sets $E$, for which introduce some simpler notation: if $E=\{d\in \N: d|n\}$ for some $n\in \N$ then $W_E=W_n$ is the ring of length $n$ Witt vectors, when $E=\{1, p, p^2, \ldots\}$ for a prime $p$ then $W_E=W^{(p)}$ is the ring of $p$-typical Witt vectors, and if $E=\N$ then $W_E=W$ is the ring of big Witt vectors.

\subsubsection{Frobenius, verschiebung, ghosts and Teichm\"uller}  Fix an index set $E$. If $p\in E$ then the Frobenius map $F_p$ is a homomorphism $W_{E}\to W_{E|p}$ and the first ghost component is a homomorphism $g: W_E\to \Ga$. Both of these are faithfully flat (as morphisms of affine schemes) and we denote by $W_E[F_p]$ and $VW_E$ their kernels. Thus, we have short exact sequences of flat affine group schemes: \[0\to VW_E \to W_E\stackrel{g}{\to} \Ga\to 0 \quad \text{ and }\quad  0\to W_E[F_p]\to W_E\stackrel{F_p}{\to} W_{E|p}\to 0.\] The kernel of the first ghost component $VW_E$ is called the verschiebung ideal.

The Frobenius maps for different $p\in E$ commute\footnote{Here I mean that the two homomorphisms $W_E\to W_{E|p\ell}$ given by $F_p\circ F_\ell$ and $F_\ell\circ F_p$ are equal.} and so we obtain Frobenius maps $F_n :W_E\to W_{E|n}$ for each $n\in E$. The $n$th ghost component for $n\in E$ is the homomorphism $g_n:=g\circ F_n$ and together these give a homomorphism \[g_E:=(g_n)_n: W_E\to \prod_{E}\Ga.\] If $S$ is a scheme over which all the primes in $E$ are invertible then the ghost map $g_E$ becomes an isomorphism when pulled back to $S$.

If $p$ is a prime then the Frobenius induces an $F_p$-linear isomorphism $VW^{(p)}\isomto (p)\otimes W^{(p)}$ whose inverse is denoted $V^{(p)}: (p)\otimes W^{(p)}\isomto VW^{(p)}$ and called the Verschiebung; we also write $V_p(-)=V^{(p)}(p\otimes-): W^{(p)}\isomto VW^{(p)}$.

Finally, the Teich\"uller map $[-]:\Ga\to W_E$ is the unique multiplicative section to $g$.

\subsubsection{The group schemes $W_E[F]$ and $W_E^\times[F]$} We define $W_E[F]=\cap_{p\in E} W_E[F_p]$ to be the intersection of the kernels of the Frobenius maps of $W_E$ and similarly $W_E^\times[F]=\cap_p W^\times_E[F_p]$. Both $W_E[F]$ and $W^\times_E[F]$ are flat affine group schemes over $\Spec(\Z)$. In fact, the restriction of the ghost map $W_E[F]\to \Ga=\Spec(\Z[t])$ identifies $W_E[F]$ with $\Spec(A)$ where $A\subset \Q[t]$ is the $\Z[t]$-algebra generated by the elements $t^{p^i}/p^{1+\cdots+p^{i-1}}$ for $p^i\in E$ with $p$ prime and the restriction of the ghost map $W_E^\times[F]\to \Gm$ induces an isomorphism $W_E^\times[F]\isomto \Spec(B)$ where $B\subset \Q[t, t^{-1}]$ is the $\Z[t, t^{-1}]$-algebra generated by the elements $(t-1)^{p^i}/p^{1+p+\cdots+p^{i-1}}$ for $p^i\in E$ with $p$ prime. For proofs of these statements see 3.3 of \cite{Drinfeld}.

An important property of $W_E[F]$ is the following: if $S=\Spec(R)$ is an affine and $\lambda \in W_E(R)$ then in $\lambda \in W_{E}[F](R)$ if and only if $\lambda\cdot VW_{E, S}=0$, and $\lambda\in VW_E(R)$ if and only if $\lambda\cdot W_{E, S}[F]=0$.

\subsubsection{Local decomposition}\label{subsec:local-decomp} Let $E$ be an index set, $S$ be an affine scheme and assume that $\ell\in E$ is invertible on $S$. Then the ghost map for $W_{E^{(\ell)}, S}$ is an isomorphism which, together with the isomorphism, $W_{E}\isomto W_{E^{(\ell)}}\circ W_{E_{(\ell)}}$ induces an isomorphism: \[W_{E, S}\isomto \prod_{n\in E^{(\ell)}} W_{E_{(\ell)}, S}.\]

Iterating this, we see that if all but a single prime $p\in E$ is invertible on $S$ then we obtain an isomorphism: \[W_{E, S}\isomto \prod_{\substack{n\in E_{(p)}}} W_{E^{(p)}, S}=\prod_{n\in E_{(p)}} W^{(p)}_S\] is an isomorphism.

In terms of this isomorphism we have the following descriptions of $W_{E, S}[F]$ and $VW_{E, S}$:

\begin{enumerate}[label=\textup{(\roman*)}]
\item  $W_{E, S}[F]=W_S^{(p)}[F_p]$, included as an additive subgroup along the factor corresponding to $1\in E$, 
\item  $W_{E, S}^\times[F]=W_S^{(p), \times}[F_p]$, included as a multiplicative subgroup along the factor at $1\in E$, and
\item \[VW_{E, S}=VW^{(p)}_S\times \prod_ {\substack{1\neq n\in E_{(p)}}} W_S^{(p)}.\]
\end{enumerate}

\subsection{The map $V(1):\mc{V}\to W$} The goal of this section is to construct an invertible $W$-module $\mc{V}$ equipped with a homomorphism of $W$-modules $V(1): \mc{V}\to W$.

So let $E$ be an index set generated multiplicatively by a finite set of primes. For $p\in E$ prime write $S_p=\Spec(\Z[1/n: n\in E,\gcd(n, p)=1])$.

Recall from (\ref{subsec:local-decomp}) that we have a decomposition of $W_{E, S_p}$ into a product of $p$-typical Witt vectors $W_{S_p}^{(p)}$: \[W_{E, S_p}\isomto \prod_{n\in E_{(p)}} W_{S_p}^{(p)}.\]

Define the invertible $W_E$-module $V_{E, S_p}$ by the formula \begin{equation} \mc{V}_{E,S_p}=((p)\otimes W_{S_p}^{(p)})\times \prod_{1\neq n\in E^{(p)}} W_{S}^{(p)}\label{eq:define-v-p}\end{equation} and define $V_{E, S_p}(1): \mc{V}_{E, S_p}\to W_{E, S_p}$ by the formula: \[V_{E, S_p}(1)=(V^{(p)}\circ (\id\otimes F_p))\times \prod_{1\neq n\in E_{(p)}} \id.\]

If $\ell\in E$ is another prime different to $p$ then it is clear from the definition (\ref{eq:define-v-p}) that the pullback of $V_{E, S_p}$ to $S_p\cap S_\ell$ is canonically isomorphic to $W_{E, S_p\cap S_\ell}$. Reversing the roles of $p$ and $\ell$ the same is true for $V_{E, S_\ell}$ and we see that under these identifications the pullbacks of the maps $V_{E, S_p}(1)$ and $V_{E, S_\ell}(1)$ are equal. As the affine schemes $S_p$ for primes $p\in E$ form an open cover of $\Spec(\Z)$ we can glue the invertible modules $\mc{V}_{E, S_p}$ and the maps $V_{E, S_p}(1)$ to obtain an a morphism of invertible $W_{E}$-modules \[V_E(1): \mc{V}_E\to W_{E}\] over $\Spec(\Z)$.

\begin{lemm} The image of $V_{E}(1)$ is contained in $VW_E$ and we have a short exact sequence \[0\to \mc{V}_E\otimes_{W_E} W_E[F]\to \mc{V}_E\stackrel{V_E(1)}{\to} VW_E\to 0.\]
\end{lemm}
\begin{proof} The claims can all be checked after base change to $S_p$ for $p\in E$, in which case it is clear by construction and the descriptions given in (\ref{subsec:local-decomp}).
\end{proof}

If $E\subset E'$ and $E'$ is also generated multiplicatively by a finite set of primes then we have a natural identification \[(V_{E'}(1): \mc{V}_{E'}\to W_{E'})\otimes_{W_{E'}} W_{E} = (V_E(1): \mc{V}_E\to W_{E}).\] Taking the limit over all such sets $E$ we obtain a morphism of $W$-modules \[V(1): \mc{V}\to W.\]

\begin{prop}\label{prop:v-is-invertible} $\mc{V}$ is an invertible $W$-module and we have a short exact sequence \[0\to \mc{V}\otimes_{W} W[F]\to \mc{V}\stackrel{V(1)}{\to} VW\to 0.\]
\end{prop}
\begin{proof} The invertibility of $\mc{V}$ is equivalent to the action of the affine group scheme $W^\times$ on the affine scheme $\underline{\Isom}_W(\mc{V}, W)$ begin free and transitive. But this can be checked after base change to $\Spec(\Z_{(p)})$ for each prime $p$ in which case it is clear from the construction of $\mc{V}$.

Exactness of the sequence can also be checked after base change to $\Spec(\Z_{(p)})$ where it follows easily using the decomposition in (\ref{subsec:local-decomp}).
\end{proof}

\begin{prop} The invertible $W$-module $\mc{V}$ is not free.
\end{prop}
\begin{proof} If $\mc{V}$ were free this would imply the existence of a Witt vector $v\in W(\Z)$ whose ghost components $g_n(v)=v_n$ satisfy $v_1=0$, $|v_{p^r}|=|p|$ for $p$ prime and $|v_n|=1$ for $n$ composite. However, we have equalities for all primes $p, \ell$: $v_{p\ell}=v_\ell\bmod p$ which is clearly impossible.
\end{proof}

\subsection{The ring scheme $\wt{W}$} We write $\wt{W}=W/W[F]$ for the sheaf of rings given by the quotient of $W$ by the ideal scheme $W[F]$. Thus, we have an exact sequence \begin{equation} 0\to W[F] \to W \to \wt{W}\to 0. \label{eq:exact-w-w1}
\end{equation}

\begin{prop}\label{coro:vw-is-invertible} $VW$ is an invertible $\wt{W}$-module and $\wt{W}$ is a flat, affine, formally smooth ring scheme.
\end{prop}
\begin{proof} The exact sequence in (\ref{prop:v-is-invertible}) is obtained from the exact sequence (\ref{eq:exact-w-w1}) by tensoring with the invertible $W$-module $\mc{V}$ and hence $VW\isomto \mc{V}\otimes_{W}\wt{W}$ is an invertible $\wt{W}$-module. The scheme $VW$ is flat, affine and formally smooth and hence the same is true for $\wt{W}$ by descent.
\end{proof}

\begin{prop}\label{prop:w*-to-w-etale} We have the following:

\begin{enumerate}[label=\textup{(\roman*)}]
\item The inclusion $\wt{W}^\times\to \wt{W}$ is flat, and hence weakly \'etale.
\item The map $W^{\times}\to \wt{W}^\times$ is faithfully flat with kernel $W^\times[F]$.
\end{enumerate}
\end{prop}
\begin{proof} As $\wt{W}^\times$ and $W^\times[F]$ are affine group schemes  both claims can be checked after base change to $S=\Spec(\Z_{(p)})$. In this case, we see from (\ref{subsec:local-decomp}) that the kernel of the homomorphism

\begin{equation} W_{S}= \prod_{\N_{(p)}}  W^{(p)}_S\stackrel{(F_p, \id, \id, \ldots)}{\longrightarrow} W_{S}=\prod_{\N_{(p)}}  W^{(p)}_S\label{eq:local-decomp-witt-p}\end{equation} is equal to $W_S[F]$ and therefore induces an isomorphism of ring schemes: \begin{equation} W_{S}\isomto \wt{W}_{S}.\label{eqn:w1-locally-isomorphic-to-w}\end{equation}

The claim in (i) now follows from the fact that $W^{\times} \to W$ is flat.

Restricting the map (\ref{eq:local-decomp-witt-p}) to groups of units it follows that $W^{\times}_S\to \wt{W}^\times_S$ is faithfully flat (as $F_p: W^{(p), \times}\to W^{(p), \times}$ is faithfully flat) and using (\ref{subsec:local-decomp}) we also see that kernel of this map is equal to $W^\times[F]$.
\end{proof}

\begin{lemm}\label{coro:lambda-invertible-w(1)} Let $S=\Spec(R)$ be an affine scheme and let $\lambda\in \wt{W}(R)$. Then $\lambda$ is invertible if and only if $\lambda:\wt{W}_S\to \wt{W}_S$ is a monomorphism.
\end{lemm}
\begin{proof} The claim that $\lambda$ is invertible can be checked after base change to $\Z_{(p)}$ in which case, as in the proof of (\ref{prop:w*-to-w-etale}), we have $W_S^{(1)}\isomto W_S$ and the claim follows from the lemma below.
\end{proof}

\begin{lemm} Let $S=\Spec(R)$ be an affine scheme and let $\lambda\in W(R)$. Then $\lambda$ is invertible if and only if $\lambda: W_S\to W_S$ is a monomorphism.
\end{lemm}
\begin{proof} One direction is clear. For the other, assume that $\lambda: W_S\to W_S$ is a monomorphism. The locus in $S$ where $\lambda$ is invertible given by the pullback of $W^\times\to W$ along the map $S\to W$ corresponding to $\lambda$ and is therefore a flat affine subscheme of $S$. To show that it is equal to $S$ we may assume that $S=\Spec(k)$ for an algebraically closed field $k$.

If $k$ has characteristic zero then $W_S\isomto \prod_{\N}\Garel[S]$ and the claim is clear. If $k$ has characteristic $p$ then \[W_S\isomto \prod_{\substack {n\in \N \\ (n, p)=1}} W_{S}^{(p)}\] and we are reduced to the analogous claim for $W_{S}^{(p)}$. In this case, an element $\lambda$ of $W^{(p)}(k)$ is invertible if and only if its image $\lambda_0\in k$ along the first ghost map is non-zero. Since $\lambda=[\lambda_0]+pw$ with $w\in W^{(p)}(k)$ we see that if $\lambda_0=0$ then $\lambda=pw: W_S^{(p)}\to W_S^{(p)}$ cannot be a monomorphism as, for example, the element $[\epsilon]\in W^{(p)}(k[\epsilon]/\epsilon^2)$ is killed by $p$. Therefore, $\lambda_0\neq 0$ and $\lambda$ is invertible.
\end{proof}

\section{The stack $\Sigma_{\Z}$}\label{sec:sigma}  In \cite{Drinfeld} Drinfeld constructs a $p$-adic formal stack $\Sigma_p$ (denoted $\Sigma$ in \cite{Drinfeld}). In this section we construct a stack $\Sigma_\Z$ which recovers Drinfeld's stack after $p$-adic completion.

\subsection{Hodge-Tate elements and $\Delta_\Z$} 

\subsubsection{Definition} Let $S=\Spec(R)$ be an affine scheme. We say that an element $v\in W(R)$ is Hodge-Tate if $\ker_S(v)=W_S[F]$. We denote by $W_{HT}\subset W$ the subsheaf of $W$ consisting of Hodge-Tate elements.

In terms of the decomposition in (\ref{subsec:local-decomp}) if $S$ is a $\Z_{(p)}$-scheme we have \[W_{HT, S}\isomto V^{(p)}((p)\otimes W_{S}^{(p), \times})\times \prod_{1\neq \N_{(p)}}W_{S}^{(p), \times}.\] Note that this cannot be used to deduce the representability of the sheaf $W_{HT}$ as the schemes $\Spec(\Z_{(p)})$ for primes $p$ do not form an fpqc cover of $\Spec(\Z)$. However, this does show that $W_{HT}\subset VW.$

The following gives some alternative characterisations of Hodge-Tate elements:

\begin{prop}\label{prop:hodge-tate-equivalences} Let $S=\Spec(R)$ be an affine scheme and let $v\in W(R)$. The following are equivalent:
\begin{enumerate}[label=\textup{(\roman*)}]
\item There exists a $W_S$-linear isomorphism $\alpha: \wt{W}_S\isomto VW_S$ such that $\alpha(1)=v,$
\item $\im_S(v)=VW_S$,
\item $\ker_S(v)=W_S[F]$.
\end{enumerate}
\end{prop}
\begin{proof} (i) implies (ii): Clear.

(ii) implies (iii): We must have $v\in VW(R)$ and therefore $W_S[F]\subset \ker_S(v)$. Thus, $v$ induces a $W_S$-linear epimorphism $W_S/W_S[F]=\wt{W}_S\to VW_S$ which, as $VW_S$ is an invertible $\wt{W}_S$-module, must be an isomorphism. Therefore, $\ker_S(v)=W_S[F]$.

(iii) implies (i): As noted above, using (\ref{subsec:local-decomp}) we see that $v\in VW(R)$. Therefore, $v$ induces a $W_S$-linear monomorphism $v: W_S/W_S[F]=\wt{W}\to VW_S$ which as $VW_S$ is an invertible $\wt{W}_S$-module must be an isomorphism by (\ref{coro:lambda-invertible-w(1)}).
\end{proof}

\begin{coro} We have $W_{HT}\isomto \underline{\Isom}_{\wt{W}}(\wt{W}, VW).$
\end{coro}
\begin{proof} This is (i) of (\ref{prop:hodge-tate-equivalences}). 
\end{proof}

\begin{coro}\label{coro:wht-algebraic} The inclusion $W_{HT}\to W$ factors through $VW$ and the induced map $W_{HT}\to VW$ is a flat and affine. In particular, $W_{HT}$ is a flat affine scheme over $\Spec(\Z)$.
\end{coro}
\begin{proof} By descent it is enough to prove the claims after base change to \[S=\underline{\Isom}_{\wt{W}}(\wt{W}, VW)=W_{HT}.\] It then follows from (i) of (\ref{prop:hodge-tate-equivalences}) that the map $W_{HT, S}\to VW_S$ is isomorphic to the flat affine map $\wt{W}^\times_S\to \wt{W}_S$ which gives the claim.
\end{proof}

\subsubsection{The stack $\Delta_\Z$} It is clear that $W_{HT}\subset W$ is stable under the action of scaling by $W^\times$ and we define the stack $\Delta_{\Z}:=W_{HT}/W^\times\subset W/W^\times$. Recall that if $S=\Spec(R)$ is an affine scheme then the $S$-points of $W/W^\times$ are identified with the groupoid whose objects are $(\xi: \mc{P}\to W_S)$ where $\mc{P}_S$ is an invertible $W_S$-module and $\xi$ is a $W_S$-linear map. The following two lemmas follow easily from (\ref{prop:hodge-tate-equivalences}) and (\ref{coro:wht-algebraic}).

\begin{lemm}\label{prop:deltaz-characterisation} Let $S$ be an affine scheme $(\xi: \mc{P}\to W_S)\in W/W^\times(S)$. Then the following are equivalent:
\begin{enumerate}[label=\textup{(\roman*)}]
\item $(\xi: \mc{P}\to W_S)\in \Delta_\Z(S)\subset W/W^\times(S)$,
\item $\im_S(\xi)=VW_S$,
\item $\ker_S(\xi)=\mc{P}\otimes_{W}W[F]$.
\end{enumerate}
In particular, whenever the conditions above hold, we have a short exact sequence \[0\to \mc{P}\otimes_{W_S} W_S[F]\to \mc{P}\to VW_S\to 0.\]
\end{lemm}

\begin{lemm} The map $\Delta_\Z\to W/W^\times\times_{\A^1/\Gm}B\Gm$ is a flat, affine monomorphism.\footnote{The map $W/W^\times \to \A^1/\Gm$ is induced by the first ghost component $g: W\to \Ga$.}
\end{lemm}

\subsubsection{The isomorphism $\Delta_\Z\isomto BW[F]$} The invertible $\wt{W}$-module $VW$ defines a morphism $\Spec(\Z) \to B \wt{W}^\times$ and we identify the fibre product \[BW^\times\times_{B\wt{W}^\times}\Spec(\Z)\] with the stack whose $S$-points are pairs $(\mc{P}, \rho)$ where $\mc{P}$ is an invertible $W_S$-module and $\rho$ is a $\wt{W}_S$-isomorphism $\mc{P}\otimes_{W_S}\wt{W}_S\isomto VW$. We note that as $W^{\times}\to \wt{W}^\times$ is faithfully flat with kernel $W^\times[F]$ it follows that $BW^\times\times_{B\wt{W}^\times}\Spec(\Z)$ is a $BW^\times[F]$-torsor.

If $(\xi: \mc{P}\to W_S)\in \Delta_\Z(S)$ then by (\ref{prop:deltaz-characterisation}) $\xi$ induces an isomorphism $\mc{P}\otimes_{W_S}\wt{W}_S\isomto VW_S$ and so we obtain a map \begin{equation}\Delta_\Z\to BW^\times\times_{B\wt{W}^\times}\Spec(\Z) \label{eq:tors-map}\end{equation}

\begin{prop}\label{prop:alt-description} The map $\Delta_{\Z}\to B W^\times\times_{B \wt{W}^\times} \Spec(\Z)$ is an isomorphism. In particular, $\Delta_{\Z}$ is a $BW^\times[F]$-torsor and the map $\Spec(\Z)\to \Delta_{\Z}$ corresponding to $(V(1): \mc{V}\to W)$ induces an isomorphism \[BW^\times[F]\isomto \Delta_\Z.\]
\end{prop}
\begin{proof} For an affine scheme $S$, the inverse of (\ref{eq:tors-map}) sends a pair $(\mc{P}, \rho)\in B W^\times\times_{B \wt{W}^\times} \Spec(\Z)(S)$ to $(\xi: \mc{P}\to W_S)\in \Delta_\Z(S)$ where $\xi$ is the composition \[\mc{P}\stackrel{x\mto \rho(x\otimes 1)}{\to} VW_S\to W_S.\]
The other claims are clear.
\end{proof}

\begin{rema} One can show that $W^\times[F](\Z)=1$ and so it follows that the groupoid $\Delta_\Z(\Z)$ is discrete. I do not know whether $(V(1): \mc{V}\to W)\in \Delta_\Z(\Z)$ is its only object. 
\end{rema}
 
\subsection{The stack $\Sigma_\Z$}

\subsubsection{Distinguished elements} For a ring $R$, an element $\xi\in W(R)$ is called distinguished if $\xi=[x]+v$ where $x\in R$ is nilpotent and $v\in W_{HT}(R)$. We write $W_{\dist}\subset W$ for the subsheaf of distinguished elements.

The sum of the Teichm\"uller map $[-]: \A^1\to W$ and the inclusion $VW\to W$ induce an isomorphism \[\A^1\times VW\isomto W\] and in terms of this description of $W$ we see that \[W_{\dist}\isomto \widehat{\A}^1\times W_{HT}\subset \A^1\times VW\isomto W.\]

In particular, it follows that the map $W_\dist\to W\times_{\A^1}\widehat{\A}^1$ is a flat affine monomorphism. Moreover, $W_\dist$ is stable under the action of $W^\times$ (this is proven easily using the fact that $W_{HT}\to VW$ is weakly \'etale).

\subsubsection{Definition of $\Sigma_\Z$} We define the stack $\Sigma_{\Z}:=W_\dist/W^\times\subset W/W^\times$. We have a cartesian diagram \[\xymatrix{\Sigma_\Z\ar[r] & W/W^\times\times_{\A^1/\Gm}\wh{\A}^1/\Gm\\
\Delta_\Z\ar[r]\ar[u] & W/W^\times\times_{\A^1\Gm}B\Gm\ar[u]}\] where the horizontal maps are flat, affine monomorphisms.

\begin{coro} The map $\Sigma_\Z\to \widehat{\A}^1/\Gm$ is flat and algebraic.
\end{coro}

\subsubsection{Local description of $\Sigma_\Z$} Our aim in this section to identify the restrictions of $\Sigma_\Z$ along various subschemes of $\Spec(\Z)$.

Write $W_\dist^{(p)}\subset W^{(p)}$ for the subsheaf whose value on a ring $R$ consists of those Witt vectors of the form $[x]+V_p(w)$ with $w\in W^{(p), \times}(R)$ and $x\in R$ nilpotent. It is an easy exercise to show that if $S=\Spec(R)$ is an affine $\Z_{(p)}$-scheme then under the identification \textup{(\ref{subsec:local-decomp})} we have \begin{equation} W_{\dist, S}=W_{\dist, S}^{(p)}\times\prod_{\substack{(n, p)=1\\ n\neq 1}} W^{(p), \times}_{S}.\label{eqn:local-dist}\end{equation}

\begin{prop}\label{prop:sigma-local} We have the following local descriptions of $\Sigma_{\Z}$:
\begin{enumerate}[label=\textup{(\roman*)}]
\item The natural map \[\Sigma_{\Z}\times \Spec(\Z_{(p)})\to W^{(p)}/W^{(p), \times}\times \Spec(\Z_{(p)})\] is fully faithful and induces an isomorphism \[\Sigma_{\Z}\times \Spec(\Z_{(p)})\isomto W_{\dist}^{(p)}/W^{(p), \times}\times \Spec(\Z_{(p)}).\]
\item 
The natural map \[\Sigma_{\Z}\times \Spec(\Q)\to \widehat{\A}^1/\Gm\times \Spec(\Q)\] is an equivalence.
\end{enumerate}
\end{prop}
\begin{proof} (i) This follows from the description (\ref{eqn:local-dist}).

(ii) We have an isomorphism over $\Spec(\Q)$: \[W\times {\Spec(\Q)}\isomto \prod_{n\in \N}\Ga\times \Spec(\Q).\] Under this identification $W_{\dist}\times \Spec(\Q)$ identifies with \[\widehat{\A}^1\times \prod_{\substack{n\in \N\\ n\neq 1}}\Gm \times \Spec(\Q)\] and $W^\times\times\Spec(\Q)$ with $\prod_{n\in \N} \Gm\times \Spec(\Q)$ and so the claim is clear.
\end{proof}

\begin{coro}\label{coro:sigma-local-drinfeld} We have \[\Sigma_\Z\times \Spec(\Q)\isomto \Spec(\Q)^{\dR}\] and \[\Sigma_\Z\times \Spf(\Z_p)\isomto \Sigma_p\] where $\Sigma_p$ is the stack constructed by Drinfeld in \cite{Drinfeld}.
\end{coro}
\begin{proof} The first isomorphism is just (ii) of (\ref{prop:sigma-local}). The second follows immediately from (i) of (\ref{prop:sigma-local}) and the definition of $\Sigma_p$ in 4.2.2 of \cite{Drinfeld} (there $\Sigma_p$ is denoted $\Sigma$).
\end{proof}

\section{Prismatization over $\Z$}\label{sec:prismatization} In this section we construction the prismatization functor: \[\Sch_\Z\to \St_{\Sigma_\Z}: X\mto X^{\Prism_\Z}\] and study its properties. Of course, the definition here is completely parallel to that in the $p$-adic case considered in by Bhatt--Lurie in \cite{Bhatt-Lurie} and many of the results we prove there are directly inspired by those of Bhatt--Lurie.

\subsection{The ring stacks $\overline{W}_E$} We denote by $\xi_{\Sigma_\Z}:\mc{P}_{\Sigma_\Z}\to W_{\Sigma_\Z}$ the pull-back of the universal invertible $W$-module equipped with a map to $W$ coming from the map $\Sigma_\Z\to W/W^\times$. For an index set $E$, write \[(\xi_{\Sigma_\Z}: \mc{P}_{\Sigma_\Z}\to W_{\Sigma_\Z})\otimes_{W_{\Sigma_\Z}} W_{E, \Sigma_\Z}=(\xi_{E, \Sigma_\Z}:\mc{P}_{E, \Sigma_\Z}\to W_{E, \Sigma_\Z}).\]

Thus $\mc{P}_{E, \Sigma_\Z}$ is an invertible $W_{E, \Sigma_\Z}$-module and $\xi_{E, \Sigma_\Z}:\mc{P}_{E, \Sigma_\Z}\to W_{E, \Sigma_\Z}$ is a quasi-ideal. Thus we may define the ring stack $\overline{W}_{E}:=\Cone(\xi_{E, \Sigma_\Z})$. If $\Spec(R)\to \Sigma_\Z$ is a morphism we just write $\overline{R}$ for $\overline{W}_1(R)$.

\begin{prop}\label{prop:w-etale} Let $S=\Spec(R)\to \Sigma_\Z$ be a map.
\begin{enumerate}[label=\textup{(\roman*)}] 
\item The natural map \[\Cone(\mc{P}_{E, \Sigma_\Z}(R)\to W_E(R))\to \overline{W}_E(R)\] is an isomorphism.
\item For $n\in \N$ and $R\to R'$ ind-\'etale the natural map \[\overline{W}_n(R)\otimes_{W_n(R)}^{\mathbf{L}} W_n(R')\to \overline{W}_n(R')\] is an isomorphism.
\item For $n\in \N$ the natural maps \[\overline{W}_n(R) \to  \overline{R} \leftarrow R\] are surjective on $\pi_0$ with nilpotent kernels.
\item The maps in \textup{(iii)} induce equivalences of categories of \'etale algebras \[\mathrm{Et}_{\overline{W}_n(R)}\leftarrow \mathrm{Et}_{\overline{R}} \to \mathrm{Et}_R.\]
\end{enumerate}
\end{prop}
\begin{proof} (i) This follows from the fact that if $\mc{P}$ is an invertible $W_{E, S}$-module and $S$ is affine then $H^1(S, \mc{P})=0$.

(ii) This follows from (i) and the fact that $R\mto W_n(R)$ commutes with \'etale base change and filtered colimits (see \cite{borger-witt-i}).

(iii) For the second claim, it is enough to show that the (obviously surjective) map \[\pi_0(\overline{W}_n(R))\to \pi_0(\overline{R})\] has nilpotent kernel.

By (ii), we may first replace $S=\Spec(R)$ with an ind-\'etale cover and assume that $\mc{P}\isomto W_S$ is free, so that $\xi=[x]+w$ where $w\in W_{HT}(R).$ Thus, we are reduced to showing that the kernel of \[W_n(R)/([x]+w)\to R/(x)\] is nilpotent. As $x$ is nilpotent, and $W_n(R)\to W_n(R/x)$ has nilpotent kernel, we may assume that $x=0$. Thus we need to show that $VW_n(R)/(w)$ is nilpotent. However, every element of $VW_n(R)$ is of the form $w\cdot a$ for some $a\in \wt{W}(R)$ and we see that $(w \cdot a)(w\cdot b)=w(w\cdot ab)$ with $(w\cdot ab)\in W(R)$ and hence $VW_n(R)^2\subset (w)$.
\end{proof}

\subsubsection{Local description of $\overline{W}_E$} Let $E$ be as in (\ref{subsec:e-typical-witt-vectors}). Let $P$ be a set of primes and write $E_{(P)}=\{n\in E: (n, p)=1 \textup{ for all }p\in P\}$ and let $S=\Spec(\Z[1/p: p\in P])$.

\begin{prop}\label{prop:local-w-ring-stacks} We have the following local descriptions of $\overline{W}$.
\begin{enumerate}[label=\textup{(\roman*)}]
\item The natural map \[\overline{W} \times \Spec(\Z_{(p)})\isomto \overline{W}^{(p)}\times \Spec(\Z_{(p)})\] is an isomorphism.
\item The natural map \[\overline{W}\times \Spec(\Q)\isomto \overline{W}_1\times \Spec(\Q)\] is an isomorphism.
\end{enumerate}
\end{prop}
\begin{proof} In both cases, as these ring stacks are flat and algebraic over $\Sigma_\Z$, to check that the given maps are isomorphisms it is enough to do so after base change to $\Spec(\Z)\to\Delta_\Z\subset \Sigma_\Z$.

(i) Let $S=\Spec(\Z_{(p)})$. By the remarks above we are reduced to showing that the natural map \[\Cone(V(1):\mc{V}_{S}\to W_{S})\to \Cone(V^{(p)}(1):\mc{V}_S^{(p)}\to W_S^{(p)}):=\Cone(V(1):\mc{V}_{S}\to W_{S})\otimes_{W}W^{(p)}\] is an isomorphism.

In this case, we have \[W_{S}\isomto \prod_{\N^{(p)}} W_S^{(p)} \quad \text{ and }\quad \mc{V}^{(p)}_{S}\isomto (p)\otimes W^{(p)}_{S}\times\prod_{\N^{(p)}-\{1\}} W_{S}^{(p)}\] with $V^{(p)}(1)=(V^{(p)}(1)_n)_{n\in \N_{(p)}}$ given by $V^{(p)}(1)_n=1$ for $n\neq 1$ and $V^{(p)}(1)_1=V^{(p)}\circ (\id\otimes F_p)$. The claim then follows as the components of $V^{(1)}(1)_n$, $n\in \N_{(p)}$ are invertible for all $n\neq 1.$

(ii) The argument is completely analogous to (ii).
\end{proof}

\begin{coro}\label{coro:comparisons} Under the identifications of \textup{(\ref{coro:sigma-local-drinfeld})}, we have \[\overline{W}\times \Spec(\Q)\isomto \Ga^{\dR} \quad \text{ and } \quad \overline{W}\times \Spf(\Z_p)\isomto \mc{R}_{\Sigma_p}\] where $\Ga^{\Prism_p}$ is the ring stack constructed by Drinfeld in \cite{Drinfeld}.
\end{coro}
\begin{proof} The first claim is clear from (ii) of (\ref{prop:local-w-ring-stacks}). The second claim follows from (i) of (\ref{prop:local-w-ring-stacks}) combined and the definition of $\mc{R}_{\Sigma_p}$ given in 1.4.1 of \cite{Drinfeld} (there denoted $\mc{R}_{\Sigma}$).
\end{proof}
 
\subsection{Definition and algebraicity of $X^{\Prism_\Z}$} We define the prismatization of a scheme $f: X\to \Spec(\Z)$ to be the prestack over $\Sigma_{\Z}$ whose $\Spec(R)\to \Sigma_\Z$-points are given by: \[\Hom_{\Sigma_\Z}(\Spec(R), X^{\Prism_\Z}) := \lim_n X(\overline{W}_n(R)).\]
It is clear that $\Spec(\Z)^{\Prism_\Z}=\Sigma_\Z$ and that $\mathbf{G}_{a}^{\Prism_\Z}=\overline{W}$. In particular, $\mathbf{G}_{a}^{\Prism_\Z}$ is algebraic over $\Sigma_\Z$: indeed, the map \[W_{\Sigma_\Z}\to (\xi_{\Sigma_\Z}:\mc{P}_{\Sigma_\Z}\to W_{\Sigma_\Z})=\mathbf{G}_{a}^{\Prism_\Z}\] is affine and faithfully flat.

\begin{prop}\label{prop:basic-props-sigma} The functor $X\mto X^{\Prism}_\Z$ commutes with Tor-independent limits, sends \'etale maps and covers to schematic \'etale maps and covers (likewise for open immersions and open covers). Moreover, $X^{\Prism_\Z}$ is a stack for the $\textup{fpqc}$ topology.
\end{prop}
\begin{proof} That the functor commutes with Tor-independent limits is clear from the definition.

Using (\ref{lemm:schematic}) we reduce the claim about \'etale morphisms and open immersion to the case where both $X=\Spec(B)$ and $Y=\Spec(A)$ are affine.  Let $\Spec(R)\to \Spec(A)^{\Prism_\Z}$ be a morphism corresponding to $h:A\to \overline{W}(R)=\lim_n \overline{W}_n(R)$. The claim boils down to the representability of the functor on $R$-algebras \[R'\mto \lim_n \Hom_A(B, \overline{W}_n(R'))\] by an \'etale $R$-algebra, where $\overline{W}_n(R')$ is viewed as an $A$-algebra via $A\to \overline{W}_{n}(R)$. But it is clear that this functor is represented by $R\to R_{B/A}$ where $R_{B/A}$ is the unique \'etale $R$-algebra corresponding to $\overline{R}\otimes_A{B}$ under the equivalences in (\ref{prop:w-etale}).

It is clear from (\ref{prop:w-etale}) that $X^{\Prism_\Z}$ is a sheaf for the \'etale topology (and hence Zariski topology), and is a stack for the fpqc topology whenever $X$ is affine or an open subscheme of an affine scheme, by the previous paragraph.

We need to show that for all faithfully flat maps of affine schemes $S\to T$ over $\Sigma_{\Z}$ that the induced map \begin{equation} X^{\Prism_\Z}(T)\to \lim X^{\Prism_\Z}(S^\bullet)\label{eq:prismatic-sheaf}\end{equation} is an isomorphism where $S^\bullet \to T$ is the Cech nerve of $S\to T.$ This can be proven exactly as in Lemma 7.3 of \cite{Bhatt-Lurie}; there the full faithfulness of the functor (\ref{eq:prismatic-sheaf}) is proven using Tannaka duality, while the essential surjectivity is rather elementary. Let us here also give an elementary proof of the full faithfulness of (\ref{eq:prismatic-sheaf}).

The full faithfulness of the functors (\ref{eq:prismatic-sheaf}) (for all $S\to T\to \Sigma$) is equivalent to the fibres of the map \[X^{\Prism_\Z}\to X^{\Prism_\Z}\times_{\Sigma_\Z} X^{\Prism_\Z}\] over morphisms from affine schemes being fpqc sheaves. However, the fibres are automatically \'etale sheaves and so for this it is enough to check over an \'etale cover of the target. Thus, choosing an affine open cover $X=\cup_{i\in I} X_i$ we are reduced to showing that the fibres of the maps \[(X_i\times_X X_j)^{\Prism_\Z}\isomto X_i^{\Prism_\Z}\times_{X^{\Prism_\Z}} X_j^{\Prism_\Z}\to X_i^{\Prism_\Z}\times_{\Sigma_\Z} X_j^{\Prism_\Z}\] are fpqc sheaves for $i, j\in I$. However, the $X_i$ are affine and the $X_i\times_X X_j$ are open subschemes of affine schemes and so it follows that both the source and the target are fpqc stacks.
\end{proof}

\begin{lemm}\label{lemm:closed-im} If $Z\to X$ is a closed immersion of schemes then $Z^{\Prism_\Z}\to X^{\Prism_\Z}$ is affine.
\end{lemm}
\begin{proof} By (\ref{prop:basic-props-sigma}) we may assume that $X=\Spec(A)$ and $Z=\Spec(A/I)$. In this case, choosing a surjection $\Z[x_j: j\in J]\to A$, we see that $\Spec(A/I)^{\Prism_\Z}\to \Spec(A)^{\Prism_\Z}$ is affine if the maps $\Spec(A/I)^{\Prism_\Z}\to (\Ga^J)^{\Prism_\Z}$ and $\Spec(A)^{\Prism_\Z}\to (\Ga^J)^{\Prism_\Z}$ are affine.

In other words, we may assume that $A=\Z[x_j: j\in J]$ is a polynomial ring. In this case, we have a flat affine epimorphism \[W_{\Sigma_\Z}^{J}\to (\Ga^J)^{\Prism_\Z}=(\Ga^{\Prism_\Z})^J\] and it is then enough to show that the morphism \begin{equation} \label{eq:closed-affine}W_{\Sigma_\Z}^J\times_{(\Ga^J)^{\Prism_\Z}}\Spec(A)^{\Prism_\Z}\to W_{\Sigma_\Z}^J\end{equation} is affine.

If $S=\Spec(R)\to W_{\Sigma_\Z}$ is a morphism corresponding to a homomorphism $h:\Z[x_j: j\in J]\to W(R)$ and a distinguished quasi-ideal $\mc{P}\to W_S$ then the $\Spec(R)$-points over $W_{\Sigma_\Z}^J$ of the fibre product in (\ref{eq:closed-affine}) is given by the groupoid of commutative diagrams \[\xymatrix{\Z[x_j: j\in J]\ar[r]^-h\ar[d] & W(R)\ar[d]\\
\Z[x_j: j\in J]/I\ar@{.>}[r] & \overline{W}(R)}\] where the vertical maps are the obvious ones.

Taking fibres of the vertical morphisms in the diagram above, we see that the groupoid of such commutative diagrams is equivalent to the groupoid of commutative diagrams \[\xymatrix{I\ar@{.>}[r]\ar[d] & \mc{P}(R) \ar[d]\\
\Z[x_j: j\in J] \ar[r]^-h & W(R)}\] where the vertical maps are the obvious ones. It follows that \[(W_{\Sigma_\Z}^J\times_{(\Ga^J)^{\Prism_\Z}}\Spec(A)^{\Prism_\Z})(R)\isomto \Hom_{\Z[x_j: j\in J]}(I, \mc{P}(R)).\] It is now clear, choosing a presentation of $I$ and using the affineness of $\mc{P}$, that the morphism (\ref{eq:closed-affine}) is affine.
\end{proof}

\begin{coro} For each scheme $X$ the morphism $X^{\Prism_\Z}\to \Sigma_\Z$ is algebraic.
\end{coro}
\begin{proof} Using (\ref{prop:basic-props-sigma}) we may assume that $X=\Spec(A)$ is affine in which case it follows from (\ref{lemm:closed-im}).
\end{proof}

\subsection{The Hodge-Tate stack $X^{HT}$} This section is directly inspired by parts of Section 5 of \cite{Bhatt-Lurie}, in particular those dealing with the Hodge-Tate stack.

The Hodge-Tate stack of a scheme $X$ is the $\Delta_\Z$-stack: \[X^{HT}:=X^{\Prism_\Z}\times_{\Sigma_\Z} \Delta_\Z.\] Thus if $\Spec(R)\to \Delta_\Z$ is a map we have \[X^{HT}(R)\isomto \lim_n X(W_{n, \Delta_\Z}(R)).\]

The (epi)morphisms of quasi-ideals $\mc{P}_{n, \Delta_\Z}\to VW_{n, \Delta_\Z}$ induce maps \[\overline{W}_{n, \Delta_\Z}=\Cone(\mc{P}_{n, \Delta_\Z}\to W_{n, \Delta_\Z})\to \Cone(VW_{n, \Delta_\Z}\to W_{n, \Delta_\Z})\isomto \Ga\] which then give maps \[\rho_{X}: X^{HT}\to X\times \Delta_\Z.\] More generally, if $f: X\to Y$ is a morphism then we denote by $\rho_{X/Y}$ the morphism \[\rho_{X/Y}: X^{HT}\to Y^{HT}\times_Y X.\]

\begin{lemm}\label{lemm:etale-localisation-ht} Let \[\xymatrix{X'\ar[r]\ar[d] & Y'\ar[d]\\
X\ar[r] & Y}\] be a commutative diagram of schemes and assume that the vertical morphisms are \'etale. Then the diagram: \[\xymatrix{X'^{HT}\ar[r]\ar[d] & Y'^{HT}\times_{Y'} X' \ar[d]\\
X^{HT} \ar[r] & Y^{HT}\times_Y X}\] is cartesian.
\end{lemm}
\begin{proof} This is straightforward to prove using the fact that $X\mto X^{HT}$ preserves \'etale morphism and that for $\Spec(R)\to \Delta_\Z$, the map \[\overline{W}_{\Delta_\Z}(R)=\Ga^{HT}(R)\to R\] has the structure of a square-zero extension (see below).
\end{proof}

Over $\Delta_\Z$ we write $W_n[F](-1):=W[F]\otimes_W \mc{P}_{n, \Delta_\Z}$ and note that $W_n[F](-1)$ is locally isomorphic to $W_{n, \Delta_\Z}[F]$ as a $W_{n, \Delta_\Z}$-module.

The homomorphism of ring stacks \[\overline{W}_{n, \Delta_\Z}\to \Ga\] admits a unique structure of a square zero extension of $\Ga$ by $B W_{n, \Delta_\Z}[F](-1)$ (following identical arguments as in \cite{Bhatt-Lurie}). This implies that, if $\Spec(R)\to Y^{HT}\times_Y X$ is a morphism corresponding to a commutative diagram \[\xymatrix{\Spec(R)\ar[r]^i\ar[d] & X\ar[d]\\
\colim_n\Spec(\overline{W}_{n, \Delta_Z}(R))\ar[r]\ar@{.>}[ur] & Y}\] then the groupoid of dotted arrows making the square commute admits a natural action of \[\lim_n \Hom_{R}(Li^*L_{X/Y}, BW_n[F](-1)(R))\isomto \Hom_{R}(Li^*L_{X/Y}, BW[F](-1)(R))\] where $L_{X/Y}$ is the cotangent complex of $X\to Y$, which is free and transitive whenever this groupoid is non-empty.

Let us write $\Theta_{X/Y}^{HT}$ for the sheaf of animated groups over $Y^{HT}\times_Y X$ given by \[\Theta_{X/Y}^{HT}(R)=\Hom_{R}(i^*(L_{X/Y}), BW[F](-1)(R)).\] We note that the compatibility of the cotangent complex with \'etale localisation shows that, in the situation of (\ref{lemm:etale-localisation-ht}), the pullback of $\Theta_{X/Y}^{HT}$ along the map $Y'^{HT}\times_{Y'}X'\to Y^{HT}\times_Y X$ is canonically isomorphic to $\Theta_{X'/Y'}^{HT}$.

The following is a reformulation of the statements above:

\begin{coro}\label{coro:pseudo-torsor} $\pi_{X/Y}: X^{HT}\to Y^{HT}\times_Y X$ is a pseudo-torsor under the sheaf of animated groups $\Theta_{X/Y}^{HT}$.
\end{coro}

\begin{prop} If $f: X\to Y$ is a syntomic morphism of schemes then:
\begin{enumerate}[label=\textup{(\roman*)}]
\item $\Theta_{X/Y}^{HT}$ is a flat algebraic stack over $Y^{HT}\times_X Y$.
\item The morphism $f^{\Prism_\Z}\to Y^{\Prism_\Z}$ is flat.
\item The morphism $\rho_{X/Y}$ is a torsor under $\Theta_{X/Y}^{HT}$ 
\end{enumerate}
\end{prop}
\begin{proof} (i) By the compatibility of $\Theta_{X/Y}^{HT}$ with \'etale localisation in $X$ and $Y$ we reduce to the case where both $X=\Spec(B)$ and $Y=\Spec(A)$ are affine. In this case $L_{B/A}\isomto [P_1\to P_0]$ is isomorphic to a two term complex whose terms are finite projective $B$-modules. This induces a cartesian square: \[\xymatrix{ P_0^{\vee}\otimes_{B} BW[F](-1) \ar[r] & P_1^{\vee}\otimes_{B} BW[F](-1)\\
\Theta_{X/Y}^{HT}\ar[r]\ar[u] & Y^{HT}\times_{Y} X\ar[u]_{0}.}\] Therefore, \[\Theta_{X/Y}^{HT}\isomto \Cone(P_0^{\vee}\otimes_B W[F](-1)\to P_1^\vee\otimes_B W[F](-1))\] and it follows that $\Theta_{X/Y}^{HT}$ is the quotient of a flat affine commutative group scheme by the action of flat affine commutative group scheme and is therefore algebraic.

(ii) We reduce immediately to the case where both $X$ and $Y$ are affine. Moreover, \'etale locally on the source and target every syntomic morphism $A\to B$ can be written in the form $A\to A[x_1, \ldots, x_n]/(f_1, \ldots, f_m)$ where $(f_1, \ldots, f_m)$ is a Koszul regular sequence. Passage to the limit then allows us to reduce to the case where $A$ is of finite type over $\Z$ and then by \cite[Tag 07M8]{stacks-project} we may assume that $A$ is a finite type polynomial algebra over $\Z$.

As $\Delta_\Z\to \Sigma_\Z$ is a nilpotent immersion, to show that $f^{\Prism_\Z}$ is flat it will be enough to show that the map $f^{HT}$ is flat, and then as $f$ is flat, that the map $\rho_{X/Y}: X^{HT}\to  Y^{HT}\times_Y X$ is flat. By (\ref{coro:pseudo-torsor}), if $\rho_{X/Y}$ is an epimorphism for the fpqc topology then it is a torsor under the algebraic group stack $\Theta_{X/Y}^{HT}$ from which it follows that $\rho_{X/Y}$ is flat. Therefore, we need only show that $\rho_{X/Y}$ is an fpqc epimorphism.

So fix a map $S=\Spec(R)\to \Sigma_\Z$ corresponding to some $(\xi: \mc{P}\to W_S)$. Then $\Hom_{\Sigma_\Z}(S, Y^{HT}\times_Y X)$ is the groupoid whose objects are commutative diagrams of the form \begin{equation} \xymatrix{A\ar[r]\ar[d] & \overline{W}(R)\ar[d]\\
B\ar[r] & R}\label{eqn:xyht-morphisms} \end{equation} while $\Hom_{\Sigma_{\Z}}(S, X^{HT})$  is the groupoid of commutative diagrams \begin{equation}\xymatrix{A\ar[r]\ar[d] & \overline{W}(R)\ar[d]\\
B\ar[r]\ar@{.>}[ur] & R.}\label{eqn:xht-morphisms}\end{equation}

Thus, to show that $X^{HT}\to X\times_Y Y^{HT}$ is an fpqc epimorphism, we need to show that after possibly replacing $R$ by a faithfully flat cover every diagram as in (\ref{eqn:xyht-morphisms}) can be filled in to a diagram as in (\ref{eqn:xht-morphisms}).

So fix a diagram as in (\ref{eqn:xyht-morphisms}). As $B=A[x_1, \ldots, x_n]/(f_1, \ldots, f_m)$ and $A$ is a polynomial algebra over $\Z$ we can lift the map $A[x_1, \ldots, x_n]\to R$ to map $A[x_1, \ldots, x_n]\to W(R)$ inducing the given map $A\to \overline{W}(R)$. This yields a commutative diagram \[\xymatrix{(f_1, \ldots, f_n) \ar[r]\ar[d] & VW(R)\ar[d]\\
A[x_1, \ldots, x_n]\ar[r]\ar[d] & W(R)\ar[d]\\
B\ar[r] & R}\]

As $\xi: \mc{P}\to VW_S$ is faithfully flat, after replacing $R$ with a faithfully flat cover, we may assume that the images of the elements $f_1,\ldots, f_m$ in $VW(R)$ are of the form $\xi(g_1), \ldots, \xi(g_m)$ for some $g_i\in \mc{P}(R)$. In this case, we can define a map $(f_1, \ldots, f_m)\to \mc{P}(R)$ by sending $f_i$ to $g_i$, noting that this is well defined as $(f_1, \ldots, f_m)$ is Koszul regular and $\mc{P}(R)\to W(R)$ is a quasi-ideal. Thus, we obtain a commutative diagram:

\[\xymatrix{(f_1, \ldots, f_n) \ar[r]\ar[d] & \mc{P}(R)\ar[d]\\
A[x_1, \ldots, x_n]\ar[r]\ar[d] & W(R)\ar[d]\\
B\ar[r] & R.}\] Taking the cones of the vertical maps in the top square induces an $A$-linear map $B=A[x_1, \ldots, x_n]/(f_1, \ldots, f_m)\to \overline{W}(R)$ making the diagram (\ref{eqn:xht-morphisms}) commute. Therefore, $\rho_{X/Y}$ is an fpqc epimorphism.

(iii) This follows from the proof of (ii).
\end{proof}

\begin{coro}\label{prop:gerbes-tangent-perfect} If $X\to Y$ is smooth quasi-compact and quasi-separated then $\rho_{X/Y}$ is a gerbe banded by the group scheme $T_{X/Y}\otimes_{\Ga} W[F](-1)$ and $R \rho_{X/Y+}(\mc{O}_{X^{HT}})$ is quasi-coherent and perfect.
\end{coro}
\begin{proof} It follows from the definition of $\Theta_{X/Y}^{HT}$ that we have \[\Theta_{X/Y}^{HT}\isomto B(T_{X/Y}\otimes_{\Ga} W[F](-1))\] which gives the first claim. For the second claim, we may assume that $f: X\to Y$ has $n$-dimensional fibres. We may also work fpqc locally on $Y^{HT}\times_Y X$ in which case the claim follows from the following: if $S=\Spec(R)$ is an affine scheme and $n\geq 0$, $g: BW_S[F]^n\to S$ is the classifying stack of $W_S[F]^n$ then $Rg_+(\mc{O}_X)$ is quasi-coherent and perfect. Quasi-coherence follows from (\ref{subsec:quasi-coherence}) and perfectness follows from Lemma 7.8 of \cite{Bhatt-Lurie}.
\end{proof}

\begin{coro}\label{coro:maps-to-gerbes-are-flat} Let $X$ be a smooth scheme over $\Z$ and $f: Y\to X^{\Prism_\Z}$ be a morphism of algbraic stacks over $\Sigma_\Z$. Then $f$ flat if and only if the composition $Y\times_{\Sigma_\Z}\Delta_\Z\to X^{HT}\to X\times \Delta_\Z$ is flat.
\end{coro}
\begin{proof} It is clear that $f$ is flat if and only if $Y\times_{\Sigma_\Z}\Delta_\Z\to X^{HT}$ is flat. However, $X^{HT}\to X\times \Delta_\Z$ is a torsor under the group stack $\Theta_{X/\Z}^{HT}$ which, as $X$ is smooth over $\Z$, is the classifying stack of a flat group scheme over $X\times \Delta_\Z$. The claim the follows from the general fact that if $G$ is flat group scheme over a stack $S$ and $T$ is a gerbe banded by $G$, then a morphism $S'\to T$ is flat and only if $S'\to S$ is flat.
\end{proof}

\subsection{Presentations of $X^{\Prism_\Z}$} Given a scheme $X$ we denote by $J(X)$ the sheaf whose $R$-points are given by $J(X)(R)=\lim_n X(W_n(R))$. In general $J(X)$ will no longer be a scheme, but if $X$ is affine then $J(X)$ is again an affine scheme. If $\Spec(R)\to \Sigma_{\Z}$ is a map the homomorphisms $W_n(R)\to \overline{W}_n(R)$ induce a map $J(X)\times \Sigma_\Z\to X^{\Prism_\Z}$.

\begin{prop} If $X$ is smooth and affine over $\Spec(\Z)$ the map $J(X)\times \Sigma_\Z\to X^{\Prism_\Z}$ is faithfully flat and affine.
\end{prop}
\begin{proof} This follows from (\ref{coro:maps-to-gerbes-are-flat}) as the map $J(X)\times \Delta_\Z\to X\times \Delta_\Z$ is flat.
\end{proof}

Let us describe \[X^{\Prism_\Z}_n=(J(X)\times \Sigma_\Z)\times_{X^{\Prism_\Z}} \cdots \times\times_{X^{\Prism_\Z}} (J(X)\times \Sigma_\Z)\] in the case that $X=\Spec(A)$ is affine.

Let $I_n=\ker(A^{\otimes_\Z n}\to A)$ and consider the affine scheme over $W$: \[J(\Spec(A^{\otimes_\Z n}[t, t^{-1}I_n]))\to J(\Spec(\Z[t]))\isomto W.\] If we let $\Gm$ act on $\Spec(A^{\otimes_\Z n}[t, t^{-1}I_n])$ by giving $t$ degree one, then $J(\Gm)=W^\times$ acts on $J(\Spec(A^{\otimes_{\Z} n}[t, t^{-1}I_n]))$ compatibly with the natural map to $W$. Therefore, we obtain a map \[J(\Spec(A^{\otimes_\Z n}[t, t^{-1}I_n]))/W^\times\to W/W^\times.\] It is straightforward to prove (just as in the de Rham case) that we have \[X^{\Prism_\Z}_n\isomto J(\Spec(A^{\otimes_\Z n}[t, t^{-1}I_n]))/W^\times\times_{W/W^\times}\Sigma_\Z.\]

Finally, if a smooth scheme $X$ admits a $\delta$-structure, which is to say a family of commuting Frobenius lifts $(\varphi_p)_{p}$, then there is a canonical morphism $X\to J(X)$. Composing this with the map above, induces a map \[X\times \Sigma_\Z\to X^{\Prism_\Z}.\]

\begin{prop} If $X\to \Spec(\Z)$ is a smooth $\delta$-scheme then the map \[X\times \Sigma_\Z\to X^{\Prism_\Z}\] is faithfully flat.
\end{prop}
\begin{proof} Faithful flatness of the map follows from (\ref{coro:maps-to-gerbes-are-flat}).
\end{proof}
If we write \[X^{\Prism_\Z}_{\delta, n}=(X\times \Sigma_\Z)\times_{X^{\Prism_\Z}} \cdots \times_{X^{\Prism_\Z}} (X\times \Sigma_\Z)\] and $X=\Spec(A)$ as above we have natural isomorphisms: \[X^{\Prism_\Z}_{\delta, n}=(J(\Spec(A^{\otimes_\Z n}[t, t^{-1}I_n]))\times_{J(\Spec(A^{\otimes_\Z n}))} \Spec(A^{\otimes_{\Z} n}))/W^\times\times_{W/W^\times} \Sigma_\Z.\]

\subsection{Comparisons and the cohomology of $X^{\Prism_\Z}$} For each prime $p$, we denote by $\Sigma_p$ the $p$-adic formal stack, and by \[\Sch_{\Spf(\Z_p)}\to \St_{\Sigma_p}:X\mapsto X^{\Prism_p}\] the $p$-adic prismatization functor, defined by Drinfeld in \cite{Drinfeld} and by Bhatt-Lurie in \cite{Bhatt-Lurie}.

\begin{theo}\label{prop:de-rham-comp} Let $X\to \Spec(\Z)$ be a scheme. Then we have:
\begin{enumerate}[label=\textup{(\roman*)}]
\item $\Sigma\times \Spec(\Q)\isomto \Spec(\Q)^{\dR}$ and $\Sigma\times \Spf(\Z)\isomto \Sigma_p$.
\item Under the equivalences in \textup{(i)} we have natural isomorphisms \[X^{\Prism_\Z}\times \Spec(\Q)\isomto X_\Q^{\dR} \quad \text{ and } \quad X^{\Prism_\Z}\times \Spf(\Z_p)\isomto X_{\Z_p}^{\Prism_p}\] where $X_{\Z_p}=X\times \Spf(\Z_p)$.
\end{enumerate}
\end{theo}
\begin{proof} This is just a restatement of (\ref{coro:comparisons}) combined with the definition of $X\mapsto X^{\dR}$ in (\ref{subsec:de-rham-stacks}) and definition of the prismatization functor in 1.4.2 of \cite{Drinfeld}.
\end{proof}

\begin{theo}\label{prop:cohom-properties} If $f: X\to Y$ is syntomic, quasi-compact and quasi-separated then $Rf_{+}^{\Prism_\Z}(\mc{O}_X)$ is quasi-coherent. It is perfect whenever $f$ is proper and smooth.
\end{theo}
\begin{proof} Quasi-coherence of $Rf_{+}^{\Prism_\Z}(\mc{O}_X)$ follows from (ii) of (\ref{coro:pseudo-torsor}) and (\ref{subsec:quasi-coherence}).

Now assume that $f$ is smooth and proper. In order to show that $Rf_{+}^{\Prism_\Z}(\mc{O}_X)$ is perfect, is enough to show that $Rf^{HT}_{+}(\mc{O}_{X^{HT}})\in D(Y^{HT})$ is perfect, since $\Delta_\Z\to \Sigma_\Z$ is a nilpotent immersion. However, by (\ref{prop:gerbes-tangent-perfect}) we know that $R\rho_{X/Y+}(\mc{O}_{X^{HT}})$ is perfect and since the projection $f_{Y^{HT}}: Y^{HT}\times_Y X\to Y^{HT}$ is smooth and proper it follows that $Rf^{HT}_{+}(\mc{O}_{X^{HT}})\isomto R f_{Y^{HT}+}(R\rho_{X/Y+}(\mc{O}_{HT}))$ is perfect.
\end{proof}

For $X/\Z$ denote by $f^{\Prism_p}_{\Z_p}: X_{\Z_p}^{\Prism_p}\to \Sigma_p$ the structure map of the $p$-adic prismatisation of $X_{\Z_p}:=X\times \Spf(\Z_p)$.

\begin{coro}\label{coro:cohomology} Let $f: X\to \Spec(\Z)$ be quasi-compact, quasi-separated and syntomic and write $i_\Q: \Spec(\Q)^{\dR}\to \Sigma$ and $i_p: \Sigma_p\to \Sigma_\Z$ for the maps induced by \textup{(i)} of \textup{(\ref{prop:cohom-properties})}. Then we have natural identifications: \[Li_{\Q}^*Rf^{\Prism_\Z}_{+}(\mc{O}_{X^{\Prism_\Z}})\isomto Rf_{\Q+}^{\dR}(\mc{O}_{X_\Q^\dR})\in D(\Spec(\Q)^{\dR})\] and \[Li_{p}^*Rf^{\Prism_\Z}_{+}(\mc{O}_{X^\Prism})\isomto Rf_{\Z_p +}^{\Prism_p}(\mc{O}_{X^{\Prism_p}_{\Z_p}})\in D(\Sigma_p).\]
\end{coro}
\begin{proof} This follows from (\ref{prop:de-rham-comp}) and (\ref{prop:cohom-properties}).
\end{proof}

\bibliography{biblio}{}
\bibliographystyle{halpha-abbrv}

\end{document}

%% file: commands.tex
\DeclareMathSymbol{A}{\mathalpha}{operators}{`A}%
\DeclareMathSymbol{B}{\mathalpha}{operators}{`B}%
\DeclareMathSymbol{C}{\mathalpha}{operators}{`C}%
\DeclareMathSymbol{D}{\mathalpha}{operators}{`D}%
\DeclareMathSymbol{E}{\mathalpha}{operators}{`E}%
\DeclareMathSymbol{F}{\mathalpha}{operators}{`F}%
\DeclareMathSymbol{G}{\mathalpha}{operators}{`G}%
\DeclareMathSymbol{H}{\mathalpha}{operators}{`H}%
\DeclareMathSymbol{I}{\mathalpha}{operators}{`I}%
\DeclareMathSymbol{J}{\mathalpha}{operators}{`J}%
\DeclareMathSymbol{K}{\mathalpha}{operators}{`K}%
\DeclareMathSymbol{L}{\mathalpha}{operators}{`L}%
\DeclareMathSymbol{M}{\mathalpha}{operators}{`M}%
\DeclareMathSymbol{N}{\mathalpha}{operators}{`N}%
\DeclareMathSymbol{O}{\mathalpha}{operators}{`O}%
\DeclareMathSymbol{P}{\mathalpha}{operators}{`P}%
\DeclareMathSymbol{Q}{\mathalpha}{operators}{`Q}%
\DeclareMathSymbol{R}{\mathalpha}{operators}{`R}%
\DeclareMathSymbol{S}{\mathalpha}{operators}{`S}%
\DeclareMathSymbol{T}{\mathalpha}{operators}{`T}%
\DeclareMathSymbol{U}{\mathalpha}{operators}{`U}%
\DeclareMathSymbol{V}{\mathalpha}{operators}{`V}%
\DeclareMathSymbol{W}{\mathalpha}{operators}{`W}%
\DeclareMathSymbol{X}{\mathalpha}{operators}{`X}%
\DeclareMathSymbol{Y}{\mathalpha}{operators}{`Y}%
\DeclareMathSymbol{Z}{\mathalpha}{operators}{`Z}%

\mathchardef\mhyphen="2D
\DeclareMathOperator*{\colim}{colim}
\usepackage{relsize}
\usepackage{bm}
\usepackage[bbgreekl]{mathbbol}
\usepackage{amsfonts}
\DeclareSymbolFontAlphabet{\mathbb}{AMSb}
\DeclareSymbolFontAlphabet{\mathbbl}{bbold}
\newcommand{\Prism}{{\mathlarger{\mathbbl{\Delta}}}}

\newcommand{\wt}[1]{\widetilde{#1}}
\newcommand{\wh}[1]{\widehat{#1}}
\newcommand{\Z}{\mathbf{Z}}

\newcommand{\N}{\mathbf{N}}
\newcommand{\aS}{\mathcal{S}}

\newcommand{\Fun}{\mathrm{Fun}}

\newcommand{\Rees}{\mathrm{Rees}}

\newcommand{\A}{\mathbf{A}}

\newcommand{\aRing}{\mathrm{aRing}}

\newcommand{\Q}{\mathbf{Q}}

\newcommand{\mc}[1]{\mathscr{#1}}
\newcommand{\mr}[1]{\mathrm{#1}}

\newcommand{\Fil}{\mr{Fil}}
\newcommand{\Ring}{\mr{Ring}}

\newcommand{\im}{\mr{im}}

\newcommand{\id}{\mr{id}}
\newcommand{\Aff}{\mr{Aff}}
\newcommand{\dAff}{\mr{dAff}}

\newcommand{\dist}{\mr{dist}}

\newcommand{\Spf}{\mr{Spf}}
\newcommand{\Spec}{\mr{Spec}}

\newcommand{\St}{\mr{St}}

\newcommand{\Sch}{\mr{Sch}}

\newcommand{\dSt}{\mr{dSt}}
\newcommand{\dPSt}{\mr{dPSt}}

\newcommand{\mto}{\mapsto}
\newcommand{\Hom}{\mr{Hom}}

\newcommand{\Cone}{\mr{Cone}}
\newcommand{\Isom}{\mr{Isom}}

\newcommand{\dR}{\mr{dR}}

\newcommand{\isomto}{\stackrel{\sim}{\longrightarrow}}
\newcommand{\Ga}{\mathbf{G}_{\mr{a}}}

\newcommand{\reason}[2][Default answer]
  
\newcommandx*\Garel[2][1,2]{\expandafter\ifx\expandafter\relax
  \detokenize{#1}\relax \mathbf{G}_{\mathrm{a}}^{#2} \else\mathbf{G}_{\mathrm{a}, #1}^{#2}\fi}

\newcommand{\Gm}{\mathbf{G}_{\mr{m}}}